\documentclass[10pt,a4paper]{article}
\usepackage{amssymb}
\usepackage{amsmath}
\usepackage{amsthm}
\usepackage[mathscr]{eucal}
\usepackage[all]{xy}
\usepackage[english]{babel}
\usepackage{comment}
\usepackage{hhline}

\xyoption{v2}
\xyoption{2cell}
\UseTwocells

\newtheorem*{theorem*}{Theorem}
\newtheorem*{MainTheorem}{Main Result}
\newtheorem{theorem}{Theorem}[section]
\newtheorem{lemma}[theorem]{Lemma}
\newtheorem{proposition}[theorem]{Proposition}
\newtheorem*{proposition*}{Proposition}
\newtheorem{corollary}[theorem]{Corollary}
\theoremstyle{definition}
\newtheorem{definition}[theorem]{Definition}
\theoremstyle{remark}
\newtheorem{remark}[theorem]{Remark}

\newdir{ >}{{}*!/-10pt/@{>}}
\newdir{ (}{{}*!/-5pt/@^{(}}
\newdir{ )}{{}*!/-5pt/@_{(}}
\def\ddar[#1]{%
\ar@{}[#1]|(.3)*{}="A" 
\ar@{}[#1]|(.7)*{}="B" 
\ar@{=>}"A";"B" }
\def\sddar[#1]{%
\ar@{}[#1]|(.4)*{}="A" 
\ar@{}[#1]|(.6)*{}="B" 
\ar@{=>}"A";"B" }

\def\Csus{\ul{C}_*}
\def\CsusGod{\ul{{\cal C}}_*}

\def\C{\mathrm{C}}

\def\cupp{\smallsmile}
\def\bb1{\mathbf{1}}
\newcommand{\scr}{\mathscr}
\def\bHom{\mathbf{Hom}}

\def\Homi{\underline{\Hom}}
\def\Hom{\mathrm{Hom}}
\def\bbZ{\mathbb{Z}}
\def\op{\mathrm{op}}
\def\ra{\rightarrow}

\def\Sm{\mathrm{Sm}}

\def\bbE{\mathbb{E}}

\def\bbQ{\mathbb{Q}}
\def\bbZtr{\bbZ_{\mathrm{tr}}}
\def\bbG{\mathbb{G}}
\def\stimes{\wedge}
\def\Corrf{c}
\def\Ntr{{\scr N}^{\tr}}
\def\Ptr{{\scr P}^{\tr}}
\def\tr{\mathrm{tr}}
\def\Homotr{\Homi_{\tr}}
\def\Nis{\mathrm{Nis}}
\def\otimestr{\otimes_{\tr}}
\def\otimesptr{\otimes^{\mathrm{pr}}_{\tr}}
\def\bbN{\mathbb{N}}

\def\colim{\mathop{\mathrm{colim}}}

\def\Cp{\mathrm{C}^+}

\def\Cm{\mathrm{C}^-}
\def\Id{\mathrm{id}}
\def\ul{\underline}
\def\ol{\overline}
\def\fksh{\mathfrak{sh}}
\def\bC{\mathbf{C}}
\def\bCb{\mathbf{C}^{\mathrm{b}}}
\def\bCp{\mathbf{C^+}}

\def\Spec{\mathrm{Spec}}
\def\fk{\mathfrak}

\def\bbA{\mathbb{A}}
\def\Kho{\mathrm{K}}

\def\Mod{\mathrm{Mod}}
\def\twocolim{2\textrm{-}\kern-0.2em\colim}
\def\D{\mathrm{D}}
\def\CH{\mathrm{CH}}

\newcommand{\mtiny}[1]{\scriptscriptstyle{#1}}

\begin{document}
\title{Levine's motivic comparison theorem revisited}
\author{Florian Ivorra}
\maketitle
\selectlanguage{english}

\begin{abstract}
\noindent For a field of characteristic zero Levine has proved in \cite[Part I, Ch. VI, 2.5.5]{MixMotLev} that the triangulated tensor categories of motives defined in \cite{MixMotVoe} and \cite{MixMotLev} are equivalent. Using some results of \cite{IvoA}, in this paper we show that the strategy of Levine's proof can also be applied on every perfect field to the categories of triangulated motives with rational coefficients or to the pseudo-abelian hulls of the integral tensor subcategories generated by motives of smooth projective schemes.
\end{abstract}

\section*{Introduction}

Using the formalism of DG categories that allows to implement relations up to homotopy and his original approach to moving lemmas, Levine has given in \cite{MixMotLev} a construction of a triangulated category of motives over a field by generators and relations. This construction was designed in such a way that one may realize triangulated motives into cohomology theories satisfying some kind of derived Bloch-Ogus axioms. In his work \cite{MixMotVoe} Voevodsky has developped a very different approach to triangulated motives relying on the properties of the Nisnevich topology and on finite correspondences, a special family of algebraic cycles on which one can perform operations of intersection theory without passing to rational equivalence.\par
For a field of characteristic zero Levine has proved in \cite[Part I, Ch. VI, 2.5.5]{MixMotLev} that the triangulated tensor categories of motives defined in \cite{MixMotVoe} and \cite{MixMotLev} are equivalent. In this paper we explain how to apply Levine's strategy of proof to extend the previous result to the categories of triangulated motives with rational coefficients or to the pseudo-abelian hulls of the integral tensor subcategories generated by motives of smooth projective schemes on every perfect field.\par
We now run briefly through the contents of this paper which is divided in two parts. For the readers' convenience we have included in part \ref{PartI} an overview of the main features of Levine's construction of triangulated motives and sketched Levine's method. All results and constructions outlined in this part are taken from \cite{MixMotLev}. Part \ref{PartII} is the main part of the paper and divided in two sections. The first one is devoted to the construction of a triangulated commutative external product in the sense of Levine \cite[Part II, Ch. I, 2.4.1]{MixMotLev}
\begin{equation*}
\mathcal{DM}^{\mathrm{eff}}_{\fksh}(F)\xrightarrow{\Upsilon^{\mathrm{eff}}_{\fksh}}DM^{\mathrm{eff}}_{-}(F)
\end{equation*}
Although this functor enjoys very nice properties listed in the beginning of section \ref{SecMotFun} it is not the right one as it is not a tensor functor and does not extend to non effective triangulated motives. However in section \ref{SecComTr} we derive from it our main result in theorem \ref{MainTheorem} which is stated below.
\begin{MainTheorem}
Assume either that
\begin{itemize}
\item{$\mathrm{char}(F)=0$ and $A=\bbZ$, or}
\item{$\mathrm{char}(F)>0$ and $A=\bbQ$.}
\end{itemize}
There exists a triangulated tensor equivalence
\begin{equation*}
{\cal DM}(F,A)\xrightarrow{\Upsilon} DM_{gm}(F,A)
\end{equation*} 
such that for a quasi-projective scheme $X$ and an integer $n$ we have an isomorphism 
\begin{equation*}
\Upsilon(\bbZ_X(n))\simeq M(X)^*(n)
\end{equation*}
where $M(X)^*$ denote the dual of the motive of $X$ in $DM_{gm}(F,A)$. In addition the equivalence $\Upsilon$ induces a triangulated tensor equivalence
\begin{equation*}
{\cal DM}(F)^{\mathrm{pr}}\xrightarrow{\Upsilon} DM_{gm}(F)^{\mathrm{pr}}
\end{equation*}
between the pseudo-abelian hulls of the integral triangulated tensor subcategories generated by the motives of smooth projective schemes.  
\end{MainTheorem}

\subsection*{Convention and notation}

Throughout the paper we use the notation in \cite{MixMotVoe,MixMotLev} and the DG formalism for which we refer to \cite{BonKap} or \cite[Part II, Ch. II]{MixMotLev}. Although \cite{BonKap} is a standard reference, we have chosen to keep the notation of \cite[Part II, Ch. II, 1.2.7]{MixMotLev}. Therefore if $\mathsf{A}$ is an additive category and $\dagger\in\{\mathrm{b},+,-,\varnothing\}$ is a boundedness condition, $\bC^{\dagger}(\mathsf{A})$ stands for the DG category of complexes of $\mathsf{A}$ whose objects are the $\dagger$-bounded complexes and morphisms are given by
$$\Hom_{\bCb(\mathsf{A})}(A,B)^n=\prod_i\Hom_{\mathsf{A}}(A^i,B^{i+n}) $$ 
with differential $d^nf=(-1)^nf\circ d_A-d_B\circ f$. Thus $\C^{\dagger}(\mathsf{A})=Z^0\bC^{\dagger}(\mathsf{A})$ and so the morphisms in the usual category of complexes are the closed morphisms of degree zero in $\bC^{\dagger}(\mathsf{A})$. More generally if $\mathsf{A}$ is a DG category one can define $\bCb(\mathsf{A})$ as the category of one sided twisted complexes\footnote{Recall that a one sided twisted complex is a family of objects $A^i\in\mathsf{A},i\in\bbZ$ and morphisms $q^{ij}:A^i\ra A^j$ such that only a finite number of $A^i$ are non zero, $q^{ij}=0$ for $i\geq j$ and
$$dq^{ij}-\sum_kq^{kj}\circ q^{ik}=0 $$
according to the sign convention of \cite{MixMotLev}. For an an additive category viewed as a DG category we recover the usual complexes.} \cite[\S 4]{BonKap}.\par
We let $F$ be a perfect field and $A$ be either the ring of integers $\bbZ$ or the field of rationals $\bbQ$. All schemes are supposed to be noetherian and separated $F$-schemes.  We denote by $\Sm_F$ the category of smooth quasi-projective schemes and by $\Ntr$ the category of Nisnevich sheaves with transfers on smooth quasi-projective schemes. For such a scheme $X$ and integers $n,m$, recall that Levine's motivic cohomology is given by
\begin{equation*}
H^m(X,\bbZ(n))_{\mathsf{L}}=\Hom_{\mathcal{DM}(F)}(\bbZ,\bbZ_X(n)[m])
\end{equation*}
and that similarly Voevodsky's motivic cohomology is given by
\begin{equation*}
\begin{split}
H^m(X,\bbZ(n))_{\mathsf{V}}&=\Hom_{DM_{-}(F)}(M(X),\bbZ(n)[m])\\
&=\Hom_{DM^{\mathrm{eff}}_{-}(F)}(M(X),\bbZ(n)[m])
\end{split}
\end{equation*}
where the last isomorphism holds for $n$ non negative by Voevodsky's cancellation theorem \cite{Cancel}. \par Recall that a (commutative) external product \cite[Part II, Ch. I, 2.4.1]{MixMotLev} is a relaxed (symmetric) monoidal functor $\mathsf{F}:\mathsf{A}\ra\mathsf{B}$ which means that one does not assume that the maps $\boxtimes^{\mathsf{F}}_{A,B}:\mathsf{F}(A)\otimes\mathsf{F}(B)\ra\mathsf{F}(A\otimes B)$ are isomorphisms. A square
\begin{equation*}
\xymatrix{{\mathsf{A}}\ar[r]^{\mathsf{F}_1}\ar[d]_{\mathsf{F}_2} & {\mathsf{B}}\ar[d]^{\mathsf{G}_1}\\
{\mathsf{C}}\ar[r]^{\mathsf{G}_2} & {\mathsf{D}}}
\end{equation*}
whose arrows are functors between categories is said to be $2$-commutative when there exists a natural isomorphism of functors $\mathsf{G}_2\circ\mathsf{F}_2\simeq \mathsf{G}_1\circ\mathsf{F}_1$ 

\part{}\label{PartI}

\section{Levine's method}

Suppose that $\mathrm{char}(F)=0$. The idea \cite[Part I, Ch. VI, 2.1.9]{MixMotLev} used in \cite{MixMotLev} to construct a triangulated tensor equivalence 
\begin{equation*}\label{EqFun}
{\cal DM}(F)\xrightarrow{\Psi} DM_{gm}(F) 
\end{equation*}
is to produce this functor as a special realization when one views $DM_{gm}(F)$ as a kind of Bloch-Ogus twisted duality theory. Since Voevodsky's motives are homological rather than cohomological as Bloch-Ogus theories, this approach requires some good duality properties in Voevodsky's category of triangulated motives which are provided by \cite[Theorem 4.3.7]{MixMotVoe}\footnote{Using Hironaka's resolution of singularities, in {\itshape{loc.cit.}} Voevodsky shows that $DM_{gm}(F)$ is a rigid tensor category.}. 
The equivalence $\Psi$ is obtained as an extension to $\mathcal{DM}(F)$ of the twisted dual motive functor
\begin{equation}\label{DualMotiveFun}
{\renewcommand{\arraystretch}{1.5}
\begin{tabular}{rcl} 
$\Sm_F^{\op}\times\bbZ$ & $\ra$ & $DM_{gm}(F)$\\
$X(n)$ & $\mapsto$  & $M(X)^*(n)$
\end{tabular}}
\end{equation}
and thus the following square is $2$-commutative
\begin{equation*}
\xymatrix{{{\cal DM}(F)}\ar[r]^{\Psi} & {DM_{gm}(F)}\\
{\Sm_F^{\op}\times\bbZ.}\ar[u]^{X(n)\mapsto\bbZ_X(n)}\ar`r[ru][ru]_(.3){X(n)\mapsto M(X)^*(n)}& {}}
\end{equation*}
The ideal situation would have been to have a functorial complex of Nisnevich sheaves with transfers that represent the functor (\ref{DualMotiveFun}). This does not happen, however Levine manages to carry out the construction relying on the following two results due to Voevodsky\footnote{The proof of (a) relies on cdh topology and its relations to relative algebraic cycles \cite{Bivariant}.}.
\begin{enumerate}
\item{\cite[Corollaries 4.1.8 and 4.2.7]{MixMotVoe} For a smooth quasi-projective scheme $X$ and a non negative integer $n$ one has a canonical isomorphism in $DM^{\mathrm{eff}}_{-}(F)$
\begin{equation*}
\Csus\big(z_{\mathrm{equi}}(X,\bbA^n,0)\big)=\Homi(M(X),\bbZ(n))
\end{equation*}
where $\Homi$ is the partially defined internal Hom. And so the fake twisted dual motive $\Homi(M(X),\bbZ(n)) $ may be represented by an explicit complex of Nisnevish sheaves with transfers}
\item{\cite[Corollary 4.3.6]{MixMotVoe} For a smooth quasi-projective scheme $X$ of dimension at most $n$ one has a canonical isomorphism in $DM^{\mathrm{eff}}_{gm}(F)$
\begin{equation*}
\Homi(M(X),\bbZ(n))=M(X)^*(n).
\end{equation*}}
\end{enumerate}
Now suppose that $F$ is an arbitrary perfect field. In \cite[Part I, Ch. IV, 1.4]{MixMotLev} Levine has given a simple argument\footnote{We refer the readers to \cite[Appendix B]{HuberKahn} for a short summary of this argument that uses De Jong theorem on alterations \cite{DeJong}.} that shows that the category $DM_{gm}(F,A)$ is  rigid in the following cases:
\begin{itemize}
\item{$\mathrm{char}(F)=0$ and $A=\bbZ$,}
\item{$\mathrm{char}(F)>0$ and $A=\bbQ$.}
\end{itemize}
One may also replace (b) by \cite[Corollary B.2]{HuberKahn}. Therefore the crucial point is to find a replacement for (a). In this paper instead of working with the complex of Nisnevich sheaves with transfers
\begin{equation*}
\Csus\big(z_{\mathrm{equi}}(X,\bbA^n,0)\big)
\end{equation*}
provided by equidimensional cycles, one introduces, using results of \cite{IvoA}, another representative of the fake twisted dual motive which may also be used in positive characteristic. This is done in subsection \ref{ReprFakeTwDual}

\section{Overview of Levine's triangulated motives}\label{SectionOverLev}

Levine's category of triangulated motives $\mathcal{DM}(F)$ is built from generators and relations according to the following two main guiding principles.
\begin{itemize}
\item{For a Bloch-Ogus twisted duality cohomology theory\footnote{The main properties of such a theory are $\bbA^1$-homotopy invariance, well behaved cycle class maps, K\"unneth isomorphisms, semi-purity and Gysin isomorphisms.} $\Gamma$ defined via cohomology of a complex of sheaves for a suitable Grothendieck topology the category $\cal{DM}(F)$ should admit a realization functor
\begin{equation*}
\fk{R}_{\Gamma}:\mathcal{DM}(F)\ra D_{\Gamma}
\end{equation*}
taking values in the triangulated category $D_{\Gamma}$ associated to $\Gamma$.}
\item{The motivic cohomology defined by $\mathcal{DM}(F)$ has to be isomorphic to Bloch's higher Chow groups. Imposing this condition may be thought of as a way to say that the category $\mathcal{DM}(F)$ is universal for the previous property.}
\end{itemize}
The most tricky part of the construction is to implement by generators and relations cycle class maps and K\"unneth morphisms (steps 1-3). At the end of the third step one gets a DG tensor category without unit from which one obtains a triangulated tensor category without unit by passing to the homotopy category of one sided twisted complexes (step 4). The $\bbA^1$-invariance and Gysin isomorphisms are only obtained at the very end by taking a quotient of this triangulated category by a suitable thick subcategory (step 5). 

\subsection*{Step 1: Moving lemma}

Levine's original approach to moving lemmas \cite{MixMotLev,HomConFil} consists in replacing the category $\Sm_F$ by a category $\mathcal{L}(F)$ fibered over $\Sm_F$. The objets ${\cal L}(F)$ are pairs $(X,f)$ where $X$ is a smooth quasi-projective scheme and $f:X'\ra X$ is a morphism in $\Sm_F$ having a smooth section $s$. A morphism from $(X,f)$ to $(Y,g)$ in ${\cal L}(F)$ is a morphism $p:X\ra Y$ in $\Sm_F$ such that there exists a flat morphism $q$ and a commutative square 
\begin{equation*}
\xymatrix{{Y'}\ar[d]_g\ar[r]^q & {X'}\ar[d]^f\\
{Y}\ar[r]_p & {X.}}
\end{equation*}
Let $\mathcal{Z}^n(X)$ the abelian group of codimension $n$ algebraic cycles on $X$. It is only contravariant for flat morphisms, however if $\mathcal{Z}^n(X)_f$ is defined as the subgroup formed of cycles in $\mathcal{Z}^n(X)$ having a well defined pull-back along $f$ then the pull-back of cycles gives a well defined functor
{\renewcommand{\arraystretch}{1.2}\begin{equation}\label{ElCycl}
\begin{tabular}{rcl}
$\mathcal{L}(F)^{\op}\times\bbN$ & $\ra$ & $\mathrm{Mod}(\bbZ)$\\
$X(n)_f$ & $\mapsto$ & $\mathcal{Z}^n(X)_f.$
\end{tabular}
\end{equation} }
Then one adds Tate twists and forms an additive category by defining
\begin{eqnarray*}
{\cal L}(F)^{\op}\times\bbN &\ra &{\cal A}_1^{\mathrm{eff}}(F)\\
X(n)_f & \mapsto & \bbZ_X(n)_f
\end{eqnarray*}
as the free additive category on ${\cal L}(F)^{\op}\times\bbN$ such that $\bbZ_{X\amalg Y}(n)_{f\amalg g}$ is the direct sum of $\bbZ_X(n)_f$ and $\bbZ_Y(m)_g$. This construction is described by the following universal property: for each functor $\mathsf{F}:{\cal L}(F)^{\op}\times\bbN\ra\mathsf{A}$ taking values in a additive category $\mathsf{A}$ such that 
$$\mathsf{F}((X\amalg Y)(n)_{f\amalg g})=\mathsf{F}(X(n)_f)\oplus\mathsf{F}(Y(n)_g) $$
there exists a unique functor $\mathsf{F}^{\mathrm{eff}}_1$ that makes the triangle below commute
\begin{equation*}
\xymatrix{{{\cal L}(F)^{\op}\times\bbN}\ar[r]\ar[rd]_{\mathsf{F}} & {{\cal A}^{\mathrm{eff}}_1(F)}\ar[d]^{\mathsf{F}^{\mathrm{eff}}_1}\\
{} & {\mathsf{A}.}}
\end{equation*} 
The cartesian product of schemes defines a tensor structure on ${\cal A}^{\mathrm{eff}}_1(F) $ and in the sequel we denote simply by $\times$ this tensor product.  

\subsection*{Step 2: K\"unneth formula}  

The next step is to take the universal external product
\begin{eqnarray*}
{\cal A}^{\mathrm{eff}}_1(F)\xrightarrow{{\iota^c}}{\cal A}^{\mathrm{eff}}_1(F)^{\otimes,c}
\end{eqnarray*}
which is a tensor category without unit\footnote{The map $\iota^c$ is a commutative external product and a fully faithful functor.} that enjoys the following universal property: for each tensor category $\mathsf{A}$ and each commutative external product
$\mathsf{F}:{\cal A}^{\mathrm{eff}}_1(F)\ra\mathsf{A}$ there exists a unique tensor functor $\mathsf{F}^{\otimes,c}$ that makes the triangle below commute
\begin{equation*}
\xymatrix{{{\cal A}^{\mathrm{eff}}_1(F)}\ar[r]^{\iota^c}\ar[rd]_{\mathsf{F}} & {{\cal A}^{\mathrm{eff}}_1(F)^{\otimes,c}}\ar[d]^{\mathsf{F}^{\otimes,c}}\\
{} & {\mathsf{A}.}}
\end{equation*}
For objects $\Gamma,\Delta$ in ${\cal A}^{\mathrm{eff}}_1(F)$ we denote by 
\begin{equation*}
\boxtimes^{c}_{\Gamma,\Delta}:\Gamma\otimes\Delta\ra\Gamma\times\Delta
\end{equation*}
the commutative external product which is a map in ${\cal A}^{\mathrm{eff}}_1(F)^{\otimes,c}$.

\subsubsection*{Variant}

However in practise the external product $\mathsf{F}$ is not commutative but only commutative up to homotopy and all higher homotopies. Indeed the external product $\mathsf{F}$ used to realize triangulated motives in Bloch-Ogus twisted duality cohomology theories is obtained via the \textit{Alexander-Whitney} map. To overcome the lack of true commutativity of the Alexander-Whitney map is one of the most technical points in Levine's construction of realization. For this purpose Levine has developed a variant of the universal external product \cite[Part II, Ch. III, 2.1.5]{MixMotLev} which is a DG tensor category without unit that fits in a commutative triangle 
\begin{equation*}
\xymatrix{{{\cal A}_1^{\mathrm{eff}}(F)^{\otimes,\fksh}}\ar[r]^{\fk{c}} & {{\cal A}_1^{\mathrm{eff}}(F)^{\otimes,c}}\\
{} & {{\cal A}_1^{\mathrm{eff}}(F)}\ar[lu]^{\iota^{\fksh}}\ar[u]^{\iota^c}}
\end{equation*} 
where $\fk{c}$ is a DG tensor functor and $\iota^{\fksh}$ is a (non commutative) external product and a fully faithful functor. For objects $\Gamma,\Delta$ in ${\cal A}^{\mathrm{eff}}_1(F)$ we denote by 
\begin{equation*}
\boxtimes^{\fksh}_{\Gamma,\Delta}:\Gamma\otimes\Delta\ra\Gamma\times\Delta
\end{equation*}
the external product which is a map in ${\cal A}_1^{\mathrm{eff}}(F)^{\otimes,\fksh}$. This construction has the following essential properties:
\begin{itemize}
\item{Given a tensor category $\mathsf{A}$ and a commutative external product $$\mathsf{F}:{\cal A}_1^{\mathrm{eff}}(F)\ra \Delta\mathsf{A}$$ there exists a canonical DG tensor functor $\mathrm{cc}\mathsf{F}^{\otimes,\fksh} $ that makes the following diagram commute
\begin{equation}\label{DiaKunneth}
\begin{tabular}{c}
\begin{math}\xymatrix{{} & {} & {\bCp(\mathsf{A})}\\
{{\cal A}_1^{\mathrm{eff}}(F)^{\otimes,\fksh}}\ar`u[u][rru]^{\mathrm{cc}\mathsf{F}^{\otimes,\fksh}}\ar[r]^{\mathfrak{c}} & {{\cal A}_1^{\mathrm{eff}}(F)^{\otimes,c}} & {\Cp(\mathsf{A})}\ar[u]_{\textrm{inclusion}}\\
{} & {{\cal A}_1^{\mathrm{eff}}(F)}\ar[u]^{\iota^c}\ar[r]^{\mathsf{F}} \ar[ru]^{\mathrm{cc}\mathsf{F}}\ar[lu]^{\iota^{\fksh}} & {\Delta\mathsf{A}}\ar[u]_{\mathrm{cc}}}\end{math}
\end{tabular}
\end{equation}
where $\mathrm{cc}$ denotes the usual associated cochain complex functor as defined in appendix \ref{AppendixDelta}. For the explicit description of $\mathrm{cc}\mathsf{F}^{\otimes,\fksh}$ we refer to \cite[Part II, Ch. III, 2.2]{MixMotLev}}.
\item{$\fk{c}$ is a homotopy equivalence}
\end{itemize}
In the sequel $\fk{m}$ is a symbol that stands either for $c$ or $\fksh$ and we consider the DG tensor category without unit 
\begin{equation*}
\mathcal{A}^{\mathrm{eff}\fk{m}}_2(F):=\mathcal{A}^{\mathrm{eff}}_1(F)^{\otimes,\fk{m}}.
\end{equation*}
The tensor product on ${\cal A}^{\mathrm{eff}}_2(F)$ is denoted by $\otimes$.

\subsection*{Step 3: Adding cycle class maps}

Let $\bbE$ be the homotopy one point DG tensor category without unit\footnote{This category is useful to implement up to homotopy the properties of the unit in the category of triangulated motives.} defined in \cite[Part II, Ch. II, 3.1.11]{MixMotLev} and $\fk{e}$ its canonical generator. The objects of $\bbE$ are direct sums of tensor power of $\fk{e}$ and one denotes by 
\begin{equation}\label{A2E}
{\cal A}_2^{\mathrm{eff}\fk{m}}(F)[\bbE]
\end{equation}
the coproduct of ${\cal A}_2^{\mathrm{eff}\fk{m}}(F)$ and $\bbE$. Let ${\cal A}_3^{\mathrm{eff}\fk{m}}(F)$ be the DG tensor category without unit obtained from (\ref{A2E}) by adding for each non zero element $\alpha\in{\cal Z}^n(X)_f$ a closed map of degree $2n$
\begin{equation*}
[\alpha]:\fk{e}\ra\bbZ_X(n)_f.
\end{equation*} 
For the zero cycle one sets $[0]=0$. Then ${\cal A}^{\mathrm{eff}\fk{m}}_4(F)$ is the DG tensor category without unit obtained from ${\cal A}_3^{\mathrm{eff}\fk{m}}(F)$ by adding the following maps.
\begin{enumerate}
\item{For each morphism $p:(Y,g)\ra (X,f)$ in ${\cal L}(F)$ and each non zero algebraic cycle $\alpha\in{\cal Z}^n(X)_f$ one adds a morphism of degree $2n-1$ 
\begin{equation*}
h_{X,Y,[\alpha],p}:\fk{e}\ra \bbZ_Y(n)_g
\end{equation*}
with the following differential
\begin{equation*}
dh_{X,Y,[\alpha],p}=p^*\circ[\alpha]-[p^*\alpha].
\end{equation*}}
\item{For each algebraic cycles $\alpha\in{\cal Z}^n(X)_f$ and $\beta\in{\cal Z}^m(Y)_g$ one adds the morphisms of degree $2(n+m)-1$ 
\begin{eqnarray*}
h^{l}_{X,Y,[\alpha],[\beta]} & : & \fk{e}^{\otimes 2}\ra \Gamma\times\Delta \\
h^{r}_{X,Y,[\alpha],[\beta]} & : &  \fk{e}^{\otimes 2}\ra \Gamma\times\Delta
\end{eqnarray*}
where $\Gamma=\bbZ_X(n)_f$ and $\Delta=\bbZ_Y(m)_g$ with the following differentials
\begin{eqnarray*}
dh^{l}_{X,Y,[\alpha],[\beta]}& = & \boxtimes^{\fk{m}}_{\Gamma,\Delta}\circ([\alpha]\otimes[\beta])-\boxtimes^{\fk{m}}_{\Gamma\times\Delta,\bb1}\circ([\alpha\times\beta]\otimes[1])\\
dh^{l}_{X,Y,[\alpha],[\beta]} & = & \boxtimes^{\fk{m}}_{\Gamma,\Delta}\circ([\alpha]\otimes[\beta])-\boxtimes^{\fk{m}}_{\bb1,\Gamma\times\Delta}\circ([1]\otimes[\alpha\times\beta]).
\end{eqnarray*}}
\item{For each $\alpha,\beta\in{\cal Z}^n(X)_f$ and integers $a,b$ one adds a morphism of degree $2n-1$ \begin{equation*}
h_{a,b,[\alpha],[\beta]}:\fk{e}\ra \bbZ_X(n)_f
\end{equation*}
with the following differential
\begin{equation*}
dh_{a,b,[\alpha],[\beta]}=[a\alpha+b\beta]-a[\alpha]-b[\beta].
\end{equation*}}
\end{enumerate}
The next step is to kill up to homotopy all the cohomology classes that do not come from algebraic cycles. Let ${\cal A}^{\mathrm{eff}\fk{m}}_5(F)$ be the DG tensor category without unit obtained from ${\cal A}_4^{\mathrm{eff}\fk{m}}(F)$ by successively adjoining morphisms according to the following inductive process. One starts with 
\begin{equation*}
{\cal A}_5^{\mathrm{eff}\fk{m}}(F)^{(0)}:={\cal A}_4^{\mathrm{eff}\fk{m}}(F)
\end{equation*}
and then assuming the DG tensor category without unit ${\cal A}_5^{\mathrm{eff}\fk{m}}(F)^{(r-1)}$ has been formed for $r\geq 1$, the category ${\cal A}_5^{\mathrm{eff}\fk{m}}(F)^{(r)}$ is defined as a $2$-colimit
\begin{equation*}
{\cal A}_5^{\mathrm{eff}\fk{m}}(F)^{(r)}:=\twocolim_{k}{\cal A}_5^{\mathrm{eff}\fk{m}}(F)^{(r-1,k)}
\end{equation*}
where ${\cal A}_5^{\mathrm{eff}\fk{m}}(F)^{(r,k)}$ is  obtained from ${\cal A}_5^{\mathrm{eff}\fk{m}}(F)^{(r,k-1)}$ by adjoining for each closed morphism of degree $2n-r$ 
\begin{equation*}
s:\fk{e}^{\otimes k}\ra\bbZ_X(n)_f
\end{equation*}
a morphism of degree $2n-r-1$
\begin{equation*}
h_s:\fk{e}^{\otimes k}\ra\bbZ_X(n)_f
\end{equation*}
such that $dh_s=s$. The DG tensor category without unit ${\cal A}^{\mathrm{eff}\fk{m}}_5(F)$ is finally the $2$-colimit
\begin{equation*}
{\cal A}_5^{\mathrm{eff}\fk{m}}(F):=\twocolim_{r}{\cal A}_5^{\mathrm{eff}\fk{m}}(F)^{(r)}.
\end{equation*}
Let ${\cal A}^{\mathrm{eff}\fk{m}}_{\mathrm{mot}}(F)$ denote the full DG tensor subcategory without unit of ${\cal A}_5^{\mathrm{eff}\fk{m}}(F) $ generated by $\fk{e}$ and objets of the form $\bbZ_X(n)_f$.

\subsubsection*{Variant}

Let ${\cal A}^{\mathrm{eff}\fk{m}}_{\mathrm{mot}}(F)^*$ be the full DG subcategory of ${\cal A}^{\mathrm{eff}\fk{m}}_{\mathrm{mot}}(F)$ generated by the objects $\fk{e}^{\otimes a}$, $\Gamma$ and $\fk{e}^{\otimes a}\otimes\Gamma$ for $\Gamma$ in $\mathcal{A}^{\mathrm{eff}}_1(F)$. The tensor structure on $\mathcal{A}^{\mathrm{eff}}_1(F)$ extends to a tensor structure on ${\cal A}^{\mathrm{eff}\fk{m}}_{\mathrm{mot}}(F)^*$ by \cite[Part II, Ch. I, 3.5.3]{MixMotLev} and the corresponding tensor product is denoted by $\times$. In addition \cite[Part I, Ch. I, 3.1.5]{MixMotLev} one has canonical functors
\begin{equation}\label{FuncMot}
{\cal A}^{\mathrm{eff}\fk{m}}_{\mathrm{mot}}(F)^*\xrightarrow{i_{\mathrm{mot}}}{\cal A}^{\mathrm{eff}\fk{m}}_{\mathrm{mot}}(F)\xrightarrow{r_{\mathrm{mot}}}{\cal A}^{\mathrm{eff}\fk{m}}_{\mathrm{mot}}(F)^*
\end{equation}
where $i_{\mathrm{mot}}$ is the inclusion functor and $r_{\mathrm{mot}}$ is a tensor functor such that $r_{\mathrm{mot}}\circ i_{\mathrm{mot}}=\Id$. The functor $ i_{\mathrm{mot}}$ is a commutative external product, the map
\begin{equation*}
\boxtimes_{\mathrm{mot}}:\Gamma\otimes\Delta\ra\Gamma\times\Delta
\end{equation*}
where $\Gamma,\Delta $ are objects of ${\cal A}^{\mathrm{eff}\fk{m}}_{\mathrm{mot}}(F)^*$ being an extension of $\boxtimes^{\fk{m}}_{\Gamma,\Delta}$. The map above induces a natural transformation
\begin{equation*}
\boxtimes_{\mathrm{mot}}:\Id\ra i_{\mathrm{mot}}\circ r_{\mathrm{mot}}.
\end{equation*}

\subsection*{Step 4: One sided twisted complexes}

Consider now the DG tensor category without unit of one sided twisted complexes \cite[\S 4]{BonKap}\footnote{We use the notation of \cite[Part II, Ch. II, 1.2.7]{MixMotLev}.}
\begin{equation*}
\bC^{\mathrm{b},\mathrm{eff}}_{\mathrm{mot}\fk{m}}(F):=\bCb({\cal A}^{\mathrm{eff}\fk{m}}_{\mathrm{mot}}(F))
\end{equation*}
and its homotopy category 
\begin{equation}\label{Khomoteff}
\Kho^{\mathrm{b},\mathrm{eff}}_{\mathrm{mot}\fk{m}}(F):=\mathsf{Ho}\left[\bC^{\mathrm{b},\mathrm{eff}}_{\mathrm{mot}\fk{m}}(F)\right].
\end{equation}
One may also consider the DG tensor category
\begin{equation*}
\bC^{\mathrm{b},\mathrm{eff}}_{\mathrm{mot}\fk{m}}(F)^*:=\bCb({\cal A}^{\mathrm{eff}\fk{m}}_{\mathrm{mot}}(F)^*)
\end{equation*}
and its homotopy category 
\begin{equation}\label{Khomoteffetoile}
\Kho^{\mathrm{b},\mathrm{eff}}_{\mathrm{mot}\fk{m}}(F)^*:=\mathsf{Ho}\left[\bC^{\mathrm{b},\mathrm{eff}}_{\mathrm{mot}\fk{m}}(F)^*\right].
\end{equation}
The functors (\ref{FuncMot}) extend to triangulated functors
\begin{equation*}
\Kho^{\mathrm{b},\mathrm{eff}}_{\mathrm{mot}\fk{m}}(F)^*\xrightarrow{i_{\mathrm{mot}}}\Kho^{\mathrm{b},\mathrm{eff}}_{\mathrm{mot}\fk{m}}(F)\xrightarrow{r_{\mathrm{mot}}}\Kho^{\mathrm{b},\mathrm{eff}}_{\mathrm{mot}\fk{m}}(F)^*
\end{equation*}
having the same properties.

\subsubsection*{Variant}

Let ${\cal A}^{\mathrm{eff}\fk{m}+}_{\mathrm{mot}}(F)$ denote the full DG tensor subcategory without unit of ${\cal A}^{\mathrm{eff}\fk{m}}_{\mathrm{mot}}(F) $ generated by $\fk{e}$ and the objects of the form $\bbZ_X(n)_f$ where $X$ is a smooth quasi-projective scheme of dimension at most $n$ and consider the DG tensor category without unit of one sided twisted complexes 
\begin{equation*}
\bC^{\mathrm{b},\mathrm{eff}+}_{\mathrm{mot}\fk{m}}(F):=\bCb({\cal A}^{\mathrm{eff}\fk{m}+}_{\mathrm{mot}}(F))
\end{equation*}
and its homotopy category 
\begin{equation}\label{Khomoteff+}
\Kho^{\mathrm{b},\mathrm{eff}+}_{\mathrm{mot}\fk{m}}(F):=\mathsf{Ho}\left[\bC^{\mathrm{b},\mathrm{eff}+}_{\mathrm{mot}\fk{m}}(F)\right].
\end{equation}
This variant will be useful in section \ref{SecComTr}.

\subsection*{Step 5: Relations}

Let $(X,f)$ be an object in $\mathcal{L}(F)$, $\hat{X}$ a closed subscheme of $X$ and let $j:U\hookrightarrow X$ the inclusion of the open complement. We write $j^*f$ for the map $U'=X'\times_XU\ra U$ induced by $f$. For a non negative integer $n$ one defines the motive with support by
\begin{equation*}
\bbZ_{X,\hat{X}}(n)_f:=\mathsf{Mc}\left(\bbZ_X(n)_f\xrightarrow{j^*}\bbZ_U(n)_{j^*f}\right)[-1]
\end{equation*} 
where $\mathsf{Mc}$ is the mapping cone. This motive belongs to $\bC^{\mathrm{b},\mathrm{eff}}_{\mathrm{mot}\fk{m}}(F)^*$ and for a cycle $\alpha\in\mathcal{Z}^n_{\hat{X}}(X)_f$ the maps
\begin{equation*}
[\alpha]:\fk{e}\ra\bbZ_X(n)_f[2n]\qquad h_{X,U,[\alpha],j^*}:\fk{e}\ra\bbZ_U(n)_{j^*f}[2n-1]
\end{equation*} 
define a natural cycle class map with support in $\bC^{\mathrm{b},\mathrm{eff}}_{\mathrm{mot}\fk{m}}(F)^*$
\begin{equation*}
[\alpha]_{\hat{X}}:\fk{e}\ra\bbZ_{X,\hat{X}}(n)_f[2n].
\end{equation*}
The triangulated tensor category without unit $\D^{\mathrm{b},\mathrm{eff}}_{\mathrm{mot}\fk{m}}(F) $ is formed from the triangulated homotopy category (\ref{Khomoteff}) by inverting the following morphisms:
\begin{enumerate}
\item{{\itshape{Homotopy.}} Let $p:(X,f)\ra (Y,g)$ a morphism in ${\cal L}(F)$ with $p:X\ra Y$ the immersion of a closed codimension one subscheme. Let $\hat{Y}$ be a closed subscheme of $Y$ and let $\hat{X}=X\times_Y\hat{Y}$. Assume that $\hat{X}$ is smooth and that we have a commutative square
\begin{equation*}
\xymatrix{{\hat{X}}\ar[r]^{p} & {\hat{Y}}\\
{\hat{X}\times\{0\}}\ar[r]\ar[u] & {\hat{X}\times\bbA^1.}\ar[u]_{\mathrm{iso.}}}
\end{equation*}
Then invert the morphism
\begin{equation*}
p^*:\bbZ_{Y,\hat{Y}}(n)_g\ra\bbZ_{X,\hat{X}}(n)_f.
\end{equation*}}
\item{{\itshape{Excision.}} Let $(X,f)\in{\cal L}(F)$, $\hat{X}$ be a closed subscheme of $X$ and $U$ be an open subscheme containing $\hat{X}$. Invert the map
\begin{equation*}
j^*:\bbZ_{X,\hat{X}}(n)_f\ra \bbZ_{U,\hat{X}}(n)_{j^*f}
\end{equation*}
where $j:U\hookrightarrow X$ is the open immersion.}
\item{{\itshape{K\"unneth isomorphism.}} Invert the morphism
\begin{equation*}
\boxtimes^{\fk{m}}_{\Gamma,\Delta}:\Gamma\otimes\Delta\ra\Gamma\times\Delta.
\end{equation*}}
for $\Gamma=\bbZ_X(n)_f$ and $\Delta=\bbZ_Y(m)_g$.
\item{{\itshape{Gysin isomorphism.}} Let $p:(P,g)\ra (X,f)$ be a morphism in ${\cal L}(F)$ such that $p:P\ra X$ is a smooth morphism that has a section $s:X\ra P$ which is a pure codimension $d$ closed immersion. Assume that $[X]\in{\cal Z}^d(P)_g$. For an integer $n\geq d$ invert the morphism 
\begin{equation*}
\xymatrix{{\fk{e}\otimes\bbZ_X(n-d)_f[-2d]\oplus\bbZ_{P\times P,X\times P}(n)_{g\times g\amalg\Delta}}\ar[d]_{{\left(\!\sideset{_0^\alpha}{_{\Delta^*}^{-\rho}}{\mathop{\,}}\right)}}\\{\bbZ_{P\times P,X\times P}(n)_{g\times g}\oplus\bbZ_{P,X}(n)_g}}
\end{equation*}
where $\alpha$ stands for the composition
\begin{equation*}
\xymatrix@C=2cm{{\fk{e}\otimes\bbZ_X(n-d)_f[-2d]}\ar[r]^{[X]_{X}\otimes p^*}\ar[rd]_{\alpha} & {\bbZ_{P,X}(d)_g\otimes\bbZ_P(n-d)_g}\ar[d]^{\boxtimes^{\fk{m}}}\\
{} & {\bbZ_{P\times P,X\times P}(n)_{g\times g}}}
\end{equation*}
and $\rho$ for the canonical map
\begin{equation*}
\rho_{g\times g,\Delta}:\bbZ_{P\times P,X\times P}(n)_{g\times g}\ra\bbZ_{P\times P,X\times P}(n)_{g\times g\amalg\Delta}.
\end{equation*}}
\item{{\itshape{Moving Lemma.}} Let $(X,f)$ be in ${\cal L}(F)$ and $g:Z\ra X$ a morphism in $\Sm_F$. Invert the morphism
\begin{equation*}
\rho_{f,g}:\bbZ_X(n)_{f\amalg g}\ra\bbZ_X(n)_f.
\end{equation*}
given by the map $(X,f)\ra (X,f\amalg g)$ in ${\cal L}(F)$ induced by the identity.}
\item{{\itshape{Unit.}} Invert the morphism
\begin{equation*}
[1]\otimes\Id:\fk{e}\otimes\bb1\ra\bb1\otimes\bb1.
\end{equation*}}
\end{enumerate} 
The triangulated category $\D^{\mathrm{b},\mathrm{eff}}_{\mathrm{mot}\fk{m}}(F)^*$ is formed from the triangulated homotopy category (\ref{Khomoteffetoile}) by inverting the morphisms (a), (b), (d), (e) and (f) listed above\footnote{The reason why we do not invert (c) is a technical issue linked to the existence of a unit. Indeed ${\cal A}^{\mathrm{eff}\fk{m}}_{\mathrm{mot}}(F)^*$ has a unit for $\times$ and so has $\D^{\mathrm{b},\mathrm{eff}}_{\mathrm{mot}\fk{m}}(F)^*$. The tensor equivalences (\ref{EqShowUnit}) then show that the motivic category has also a unit.}. The functors (\ref{FuncMot}) extend to triangulated functors
\begin{equation*}
\D^{\mathrm{b},\mathrm{eff}}_{\mathrm{mot}\fk{m}}(F)^*\xrightarrow{i_{\mathrm{mot}}}\D^{\mathrm{b},\mathrm{eff}}_{\mathrm{mot}\fk{m}}(F)\xrightarrow{r_{\mathrm{mot}}}\D^{\mathrm{b},\mathrm{eff}}_{\mathrm{mot}\fk{m}}(F)^*
\end{equation*}
having the same properties. Now the natural transformation
\begin{equation*}
\boxtimes_{\mathrm{mot}}:\Id\ra i_{\mathrm{mot}}\circ r_{\mathrm{mot}}.
\end{equation*}
is an isomorphism by \cite[Part I, Ch. I, 3.4.2]{MixMotLev} and so 
\begin{equation}\label{EqShowUnit}
\D^{\mathrm{b},\mathrm{eff}}_{\mathrm{mot}\fk{m}}(F)^*\xrightarrow{i_{\mathrm{mot}}}\D^{\mathrm{b},\mathrm{eff}}_{\mathrm{mot}\fk{m}}(F)\qquad\D^{\mathrm{b},\mathrm{eff}}_{\mathrm{mot}\fk{m}}(F)\xrightarrow{r_{\mathrm{mot}}}\D^{\mathrm{b},\mathrm{eff}}_{\mathrm{mot}\fk{m}}(F)^*
\end{equation}
are triangulated tensor functors quasi-inverse to each others. The triangulated tensor category $\D^{\mathrm{b},\mathrm{eff+}}_{\mathrm{mot}\fk{m}}(F) $ is formed from the triangulated homotopy category (\ref{Khomoteff+}) by inverting the morphisms (a) to (f) listed above when they do belong to the category (\ref{Khomoteff+}).

\subsection*{Step 6: Inverting the Tate motive}

By \cite[Part I, Ch. VI, 2.5.4]{MixMotLev} the symmetry isomorphism
\begin{equation*}
\bbZ(1)\otimes\bbZ(1)\ra\bbZ(1)\otimes\bbZ(1)
\end{equation*} 
is the identity in $\D^{\mathrm{b},\mathrm{eff}}_{\mathrm{mot}\fk{m}}(F)$, $\D^{\mathrm{b},\mathrm{eff}}_{\mathrm{mot}\fk{m}}(F)^*$ and $\D^{\mathrm{b},\mathrm{eff}+}_{\mathrm{mot}\fk{m}}(F)$. Therefore one can invert $\bbZ(1)$ in those triangulated tensor categories to obtain the following triangulated tensor categories
\begin{equation*}
\D^{\mathrm{b}}_{\mathrm{mot}\fk{m}}(F):=\D^{\mathrm{b},\mathrm{eff}}_{\mathrm{mot}\fk{m}}(F)\Big[\bbZ(1)^{-1}\Big]\\
\end{equation*}
\begin{equation*}
\D^{\mathrm{b}}_{\mathrm{mot}\fk{m}}(F)^*:=\D^{\mathrm{b},\mathrm{eff}}_{\mathrm{mot}\fk{m}}(F)^*\Big[\bbZ(1)^{-1}\Big]
\end{equation*}
and
\begin{equation*}
\D^{\mathrm{b}+}_{\mathrm{mot}\fk{m}}(F):=\D^{\mathrm{b},\mathrm{eff}+}_{\mathrm{mot}\fk{m}}(F)\Big[\bbZ(1)^{-1}\Big].
\end{equation*}
Finally the triangulated tensor categories ${\cal DM}_{\fk{m}}(F)$, ${\cal DM}_{\fk{m}}(F)^*$ and ${\cal DM}^{+}_{\fk{m}}(F)$ are the pseudo-abelian hulls of those categories. The fact that the category ${\cal DM}(F)$ is really the category whose construction is described in \cite[Part I, Ch. I]{MixMotLev} is proved in \cite[Part I, Ch. VI, 2.5.4]{MixMotLev}.

\subsection*{General picture}

The various categories described so far belong to a commutative diagram
\begin{equation*}
\xymatrix@C=.5cm@R=.5cm{{} & {\D^{\mathrm{b}+}_{\mathrm{mot}\fksh}(F)}\ar[rr]^{\fk{c}^{+}}\ar[d] & {} & {\D^{\mathrm{b}+}_{\mathrm{mot}}(F)}\ar[d] & {}\\
{} & {\D^{\mathrm{b}}_{\mathrm{mot}\fksh}(F)}\ar[rr]^{\fk{c}}\ar[ldd]^{r_{\mathrm{mot}}} & {} & {\D^{\mathrm{b}}_{\mathrm{mot}}(F)}\ar[ldd]^{r_{\mathrm{mot}}} & {}\\
{} & {} & {\D^{\mathrm{b}}_{\mathrm{mot}\fksh}(F)^*}\ar[lu]^{i_{\mathrm{mot}}}\ar@{.>}[rr]^{\fk{c}^*} & {} & {\D^{\mathrm{b}}_{\mathrm{mot}}(F)^*}\ar[lu]^{i_{\mathrm{mot}}}\\
{\D^{\mathrm{b}}_{\mathrm{mot}\fksh}(F)^*}\ar[rr]^{\fk{c}^*} & {} & {\D^{\mathrm{b}}_{\mathrm{mot}}(F)^*} & {} & {}}
\end{equation*}
where all the arrows are triangulated tensor functors. By construction the vertical arrows are equivalences and the fact that the functors $\fk{c}$, $\fk{c}^*$ and $\fk{c}^{+}$ are equivalences too is proved in \cite[Part I, Ch. V, 1.3.6]{MixMotLev}. So by taking the pseudo-abelian hull we see that the category of triangulated motives and its variants fit in a commutative diagram of triangulated tensor equivalences
\begin{equation*}
\xymatrix@C=.5cm@R=.5cm{{} & {{\cal DM}^{+}_{\fksh}(F)}\ar[rr]^{\fk{c}^{+}}\ar[d] & {} & {{\cal DM}^{+}(F)}\ar[d] & {}\\
{} & {{\cal DM}_{\fksh}(F)}\ar[rr]^{\fk{c}}\ar[ldd]^{r_{\mathrm{mot}}} & {} & {{\cal DM}(F)}\ar[ldd]^{r_{\mathrm{mot}}} & {}\\
{} & {} & {{\cal DM}_{\fksh}(F)^*}\ar[lu]^{i_{\mathrm{mot}}}\ar@{.>}[rr]^{\fk{c}^*} & {} & {{\cal DM}(F)^*}\ar[lu]^{i_{\mathrm{mot}}}\\
{{\cal DM}_{\fksh}(F)^*}\ar[rr]^{\fk{c}^*} & {} & {{\cal DM}(F)^*} & {} & {}}
\end{equation*}

\part{}\label{PartII}

\section{The motivic commutative external product}\label{SecMotFun}

In this section we construct a triangulated commutative external product
\begin{equation}\label{FunTrExPr}
\mathcal{DM}^{\mathrm{eff}}_{\fksh}(F)\xrightarrow{\Upsilon^{\mathrm{eff}}_{\fksh}}DM^{\mathrm{eff}}_{-}(F)
\end{equation}
which has the following properties.
\begin{enumerate}
\item{The square 
\begin{equation}\label{FunSq}
\begin{tabular}{c}
\begin{math}
\xymatrix@C=1.5cm{{\mathcal{DM}^{\mathrm{eff}}_{\fksh}(F)}\ar[r]^{\Upsilon_{\fksh}^{\mathrm{eff}}} & {DM^{\mathrm{eff}}_{-}(F)}\\
{{\cal L}(F)^{\op}\times\bbN}\ar[r]\ar[u]^{X(n)_f\mapsto\bbZ_X(n)_f} & {\Sm_F^{\op}\times\bbN.}\ar[u]_{X(n)\mapsto \Homi(M(X),\bbZ(n))}}
\end{math}
\end{tabular}
\end{equation}
is $2$-commutative.}
\item{Let $\Gamma$ be an object in $\mathcal{DM}^{\mathrm{eff}}_{\fksh}(F)$ and $n$ a non negative integer. The induced map 
\begin{equation*}
\Hom_{\mathcal{DM}^{\mathrm{eff}}_{\fksh}(F)}\big(\bbZ(n),\Gamma\big)\xrightarrow{\Upsilon^{\mathrm{eff}}_{\fksh}}\Hom_{DM^{\mathrm{eff}}_{-}(F)}\big(\bbZ(n), \Upsilon^{\mathrm{eff}}_{\fksh}(\Gamma)\big)
\end{equation*}
is an isomorphism.}
\end{enumerate}

\subsection{Representing the fake twisted dual motive}\label{ReprFakeTwDual}

We shall now explain how to represent the fake twisted dual motive by an explicit functorial complex of Nisnevich with transfers. By \cite[Proposition 3.2]{IvoA} the Godement monad induces a monad $\scr G_{\mathrm{Nis}}^{\tr}$ over the category of Nisnevich sheaves with transfers and therefore the cosimplicial Godement resolution ${\scr G}_{\Nis}^{\mathrm{tr},*}(\scr F)$ of a Nisnevich sheaf with transfers $\scr F$ has a canonical structure of cosimplicial Nisnevich sheaf with transfers. We use the notation in appendix \ref{AppendixDelta}. 

\begin{definition}\label{NotnCsusGod}
For a Nisnevich sheaf with transfers $\scr F$ we let $\CsusGod(\scr F)$ be the complex of Nisnevich sheaves with transfers
\begin{equation}
\CsusGod(\scr F)=\textsf{Tot}\Big[\Csus\mathrm{cc}{\scr G}_{\Nis}^{\mathrm{tr},*}(\scr F)\Big]
\end{equation}
where $\textsf{Tot}$ is the complex associated to a double complex as in \cite[Notation 11.5.2]{KasSch}.
\end{definition}

\noindent Given a smooth quasi-projective scheme $X$ we set for short $\CsusGod(X)=\CsusGod(\bbZtr[X])$. The usual properties of Godement resolutions give us the following remark.

\begin{remark}\label{RemFlabiness}
Given a Nisnevich sheaf with transfers ${\scr G}$ let $\Csus^{\nu}(\scr G)$ be the complex of Nisnevich sheaves with transfers
\begin{equation*}
\Csus^{\nu}(\scr G):\cdots\ra 0\ra\ul{C}_{\nu}(\scr G)\ra\ul{C}_{\nu-1}(\scr G)\ra\cdots\ra\ul{C}_1(\scr G)\ra\ul{C}_0(\scr G).
\end{equation*}
With the above notation we may write
\begin{equation*}
\CsusGod(\scr F)=\colim_{\nu}\mathsf{Tot}\Big[\Csus^{\nu}{\scr G}_{\Nis}^{\mathrm{tr},*}(\scr F)\Big]
\end{equation*}
and so $\CsusGod(\scr F)$ is a filtred colimit of bounded below complexes of flabby sheaves. 
\end{remark}

\noindent Recall that one defines \cite[\S 3]{SusVoeA} the Nisnevich sheaf with transfers $\bbZtr[\bbG_m^{\stimes n}]$ as the quotient 
$$\bbZtr[\bbG_m^{\stimes n}]:=\bbZtr[\bbG_m^{\times n}]/{\cal D}_n$$
where ${\cal D}_n$ is the subsheaf obtained by taking the union of the image of the morphisms $\bbZtr[\bbG_m^{\times n-1}]\ra\bbZtr[\bbG_m^{\times n}]$ given by the closed immersions
$$(x_1,\ldots,x_{n-1})\mapsto (x_1,\ldots,1,\ldots,x_{n-1}).$$
The image in $DM^{\mathrm{eff}}_{-}(F)$ of the Nisnevich sheaf with transfers $\bbZtr[\bbG_m^{\stimes n}]$ is the motive $\bbZ(n)[n]$.

\begin{definition}\label{fkzsh}
Let $X$ be a smooth quasi-projective scheme, $W$ a closed subscheme of $X$ and $n$ a non negative integer. We set
\begin{equation*}
\Upsilon^{\mathrm{eff}}_{W}(X(n))=\Homotr\big(\bbZtr^W[X],\CsusGod(\bbG_m^{\stimes n})\big)[-n].
\end{equation*}
where the Nisnevich sheaf with transfers $\bbZtr^W[X]$ is the cokernel of the map $\bbZtr[X\setminus W]\ra\bbZtr[X]$. 
\end{definition}

\begin{proposition}\label{lemGodRe}
Given a Nisnevich sheaf with transfers ${\scr F}$, we have a canonical isomorphism in $DM^{\mathrm{eff}}_{-}(F)$
\begin{equation*}
\Homotr(\bbZtr^W[X],\CsusGod({\scr F}))=\Homi(M_W(X),\Csus({\scr F})).
\end{equation*}
In particular we have a canonical isomorphism in $DM^{\mathrm{eff}}_{-}(F)$
\begin{equation}\label{IsolemGodRe}
\Upsilon^{\mathrm{eff}}_{W}(X(n))=\Homi(M_W(X),\bbZ(n)).
\end{equation}
\end{proposition}

\begin{proof}
By remark \ref{RemFlabiness} we have isomorphisms in  $\D(\Ntr)$
\begin{equation*}
\begin{split}
\Homotr\big(\bbZtr^W[X],\CsusGod({\scr F})\big)&=\mathsf{R}\Homotr\big(\bbZtr^W[X],\CsusGod({\scr F})\big)\\
&=\mathsf{R}\Homotr\big(\bbZtr^W[X],\Csus({\scr F})\big).
\end{split}
\end{equation*}
The isomorphism in $DM^{\mathrm{eff}}_{-}(F)$
\begin{equation*}
\begin{split}
\Homotr\big(\bbZtr^W[X],\CsusGod({\scr F})\big)&=\mathsf{R}\Homotr\big(\Csus^W(X),\Csus({\scr F})\big)\\
&=\Homi\big(M_W(X),\Csus({\scr F})\big)
\end{split}
\end{equation*}
comes from the fact that $\Csus({\scr F})$ is $\bbA^1$-local and that the morphism $\bbZtr^W[X]\ra\Csus^W(X)$ is an $\bbA^1$-quasi-isomorphism \cite[Proposition 3.2.3]{MixMotVoe}. Since  $\Csus(\bbG_m^{\stimes n})=\bbZ(n)[2n]$ we have a canonical isomorphism in $DM^{\mathrm{eff}}_{-}(F)$
\begin{equation*}
\begin{split}
\Upsilon^{\mathrm{eff}}_{W}(X(n)) & = \Homotr\big(\bbZtr[X],\CsusGod(\bbG_m^{\stimes n})\big)[-n]\\
&=\mathsf{R}\Homotr\big(\Csus(X),\Csus(\bbG_m^{\stimes n})\big)[-n]\\
&=\Homi\big(M(X),A(n)\big).
\end{split}
\end{equation*}
This proves the proposition.
\end{proof}

\noindent From the isomorphism (\ref{IsolemGodRe}) given in the previous proposition we obtain a canonical isomorphism in $\D(\bbZ)$
\begin{equation}\label{IsoCoh}
\mathsf{\Gamma}_{\mathrm{Nis}}(\Spec(F),\Upsilon^{\mathrm{eff}}_{W}(X(n)))=\mathsf{R}\mathsf{\Gamma}^{\mathrm{Nis}}_W(X,\bbZ(n)).
\end{equation}
Now by \cite[Corollary B.2]{HuberKahn} we get the following result.

\begin{corollary}\label{Dualfksh}
Let $X$ be a smooth quasi-projective scheme of dimension at most $n$ such that $M(X)$ has a dual in $DM_{gm}(F,A)$. Then the complex $\Upsilon^{\mathrm{eff}}_{W}(X(n))$ belongs to $DM^{\mathrm{eff}}_{gm}(F,A)$ and we have a canonical isomorphism
\begin{equation*}
\Upsilon^{\mathrm{eff}}_{W}(X(n))=M^*(X)(n)
\end{equation*}
 in $DM^{\mathrm{eff}}_{gm}(F,A)$. 
\end{corollary}

\noindent However to overcome the lack of true commutativity of the {\itshape{Alexander-Whitney}} map we shall not work directly with the complex given in definition \ref{fkzsh}. We consider instead the underlying cosimplicial simplicial Nisnevich sheaf with transfers:  
\begin{definition}
Let $X$ be a smooth quasi-projective scheme, $W$ a closed subscheme of $X$ and $n$ a non negative integer. We denote by $\Upsilon^{\mathrm{eff}}_{W\mathrm{cs}}(X(n))$ the object of $\mathsf{\Delta}\mathsf{\Delta}^{\op}\Ntr$
\begin{equation}\label{SimpCosimp}
\Upsilon^{\mathrm{eff}}_{W\mathrm{cs}}(X(n)):=\Homotr\big(\bbZtr^W[\Delta^*_X],\scr G_{\Nis}^{\tr,*}(\bbZtr[\bbG_m^{\stimes n}])\big)
\end{equation}
where the cosimplicial Nisnevich sheaf with transfers $\bbZtr^W[\Delta^*_X]$ is the cokernel of the map $\bbZtr[\Delta^*_X]\ra\bbZtr[\Delta^*_{X\setminus W}] $.
\end{definition}
\noindent In the cosimplicial simplicial object (\ref{SimpCosimp}) the simplicial structure is induced by the cosimplicial scheme $\Delta^*_X$ and the cosimplicial structure is induced by the cosimplicial structure of the Godement resolution. The relation between the two definitions is given by 
\begin{equation*}
\begin{split}
\mathsf{Tot}\Big[\mathrm{cc}\Big(\mathrm{c}\Upsilon^{\mathrm{eff}}_{W\mathrm{cs}}(X(n))[-n]\Big)\Big]&=\mathsf{Tot}\Big[\mathrm{cc}\Big(\Homotr(\bbZtr^W[X],\Csus{\scr H}^*_n)[-n]\Big)\Big]\\
&=\mathsf{Tot}\Big[\Homotr\Big(\bbZtr^W[X],\ul{C}_*\mathrm{cc}{\scr H}^*_n[-n]\Big)\Big]\\
&=\Homotr\Big(\bbZtr^W[X],\mathsf{Tot}\Big[\ul{C}_*\mathrm{cc}{\scr H}^*_n\Big][-n]\Big)
\end{split}
\end{equation*}
where ${\scr H}^*_n=\scr G^{\tr,*}_{\Nis}(\bbZtr[\bbG^{\stimes n}_m])$. This implies by definition that
\begin{equation*}
\mathsf{Tot}\Big[\mathrm{cc}\Big(\mathrm{c}\Upsilon^{\mathrm{eff}}_{W\mathrm{cs}}(X(n))[-n]\Big)\Big]=\Upsilon^{\mathrm{eff}}_{W}(X(n)).
\end{equation*}

\subsection{The K\"unneth formula}

\noindent We now start the construction of the triangulated functor (\ref{FunTrExPr}). For this we consider the functor which commutes with finite coproducts
\begin{equation*}
\Upsilon^{\mathrm{eff}}_{\mathrm{cs}}:\Sm_F^{\op}\times\bbN\ra\mathsf{\Delta}^{\op}\mathsf{\Delta}\Ntr.
\end{equation*}

\subsubsection{Principle of construction}\label{ParPrincipeCons}

Let ${\mathsf{A}}$ be a tensor category and assume to be given a commutative external product
\begin{equation*}
{\mathsf{F}}:{\cal L}(F)^{\op}\times\bbN\ra\mathsf{\Delta}^{\op}\mathsf{\Delta}{\mathsf{A}}
\end{equation*}
which commutes with finite coproducts. Using the shuffle map the functor $\mathrm{c}$ is a commutative external product and so by composition we obtain a commutative external product
\begin{equation*}
{\mathrm{c}}{\mathsf{F}}:{\cal L}(F)^{\op}\times\bbN\ra\mathsf{\Delta}\Cm({\mathsf{A}}).
\end{equation*}
Now consider the commutative external product obtained by a shift
\begin{eqnarray*}
\ol{\mathrm{c}\mathsf{F}}:{\cal L}(F)^{\op}\times\bbN & \ra & \mathsf{\Delta}\Cm({\mathsf{A}})\\
X(n)_f & \mapsto  & \mathrm{c}\mathsf{F}(X(n)_f)[-n].
\end{eqnarray*}
Since this functor commutes with finite coproducts, using the properties of the category ${\cal A}_1^{\mathrm{eff}}(F)$ recalled in step 1 we have a commutative triangle
\begin{equation*}
\xymatrix{{{\cal L}(F)^{\op}\times\bbN}\ar[r]\ar[rd]_{\ol{\mathrm{c}\mathsf{F}}} & {{\cal A}^{\mathrm{eff}}_1(F)}\ar[d]^{\ol{\mathrm{c}\mathsf{F}}^{\mathrm{eff}}_1}\\
{} & {\mathsf{\Delta}\Cm({\mathsf{A}})}}
\end{equation*} 
where $\ol{\mathrm{c}\mathsf{F}}^{\mathrm{eff}}_1$ is a canonical additive commutative external product. Now we may use the properties of the DG tensor category without unit ${\cal A}_2^{\mathrm{eff}\fksh}(F)$ recalled in the variant of step 2 in order to obtain a commutative diagram
\begin{equation}\label{DiaHandleKun}
\begin{tabular}{c}
\begin{math}\xymatrix{{} & {} & {\bCp(\Cm(\mathsf{A}))}\\
{{\cal A}_2^{\mathrm{eff}\fksh}(F)}\ar`u[u][rru]^{\mathrm{cc}(\ol{\mathrm{c}\mathsf{F}}^{\mathrm{eff}}_1)^{\otimes,\fksh}}\ar[r]^{\mathfrak{c}} & {{\cal A}_2^{\mathrm{eff}}(F)} & {\Cp(\Cm(\mathsf{A}))}\ar[u]_{\textrm{inclusion}}\\
{} & {{\cal A}_1^{\mathrm{eff}}(F)}\ar[u]^{\iota^c}\ar[r]^{\ol{\mathrm{c}\mathsf{F}}^{\mathrm{eff}}_1} \ar[ru]^{\mathrm{cc}(\ol{\mathrm{c}\mathsf{F}}^{\mathrm{eff}}_1)}\ar[lu]^{\iota^{\fksh}} & {\mathsf{\Delta}\Cm(\mathsf{A})}\ar[u]_{\mathrm{cc}}}\end{math}
\end{tabular}
\end{equation}
where $\mathrm{cc}(\ol{\mathrm{c}\mathsf{F}}^{\mathrm{eff}}_1)^{\otimes,\fksh}$ is a canonical DG tensor functor. The diagram (\ref{DiaHandleKun}) gives us a DG tensor functor
\begin{equation*}
\xymatrix@C=2cm{{{\cal A}_2^{\mathrm{eff}\fksh}(F)}\ar[r]_{\mathrm{cc}(\ol{\mathrm{c}\mathsf{F}}^{\mathrm{eff}}_1)^{\otimes,\fksh}}\ar@/^2em/[rr]^{\mathrm{Tot}\,\circ\,\mathrm{cc}(\ol{\mathrm{c}\mathsf{F}}^{\mathrm{eff}}_1)^{\otimes,\fksh}} & {\bCp(\Cm(\mathsf{A}))}\ar[r]_{\mathrm{Tot}} & {\bC(\mathsf{A}).}}
\end{equation*}
Recall that the right hand side $\bC(\mathsf{A})$ denotes the DG category of complexes over $\mathsf{A}$ as described at the beginning of this paper.

\subsubsection{How to handle K\"unneth morphisms}

In this subsection we apply the previous construction principle to the functor $\Upsilon^{\mathrm{eff}}_{\mathrm{cs}}$. For this we need a variant of \cite[Proposition 3.5]{IvoA} given in lemma \ref{VarTensGodtr}.
Let the ${\scr F}_n$ be a family of Nisnevich sheaves with transfers indexed by $\bbN$. A commutative external product on this family is the data of morphisms
\begin{equation*}
\boxtimes^{{\scr F_n}}_{X,Y}:{\scr F_n}(X)\otimes{\scr F_m}(Y)\ra{\scr F_{n+m}}(X\times Y)\qquad X,Y\in\mathrm{Sch}_F
\end{equation*}
which are associative commutative and functorial for finite correspondences. Denote by ${\scr Q}_{\tr}$ the category with objects the family of Nisnevich sheaves with transfers indexed by $\bbN$ endowed with a commutative external product and obvious morphisms, then the proof of \cite[Proposition 3.5]{IvoA} yields the following result.

\begin{lemma}\label{VarTensGodtr}
The monad ${\scr G}_{\mathrm{Nis}}^{\mathrm{tr}}$  induces a monad on the category $\scr Q_{\mathrm{tr}}$. Thus if the $\scr F_n$ are a family of Nisnevich sheaves with transfers endowed with a commutative external product, the ${\scr G}_{\mathrm{Nis}}^{\mathrm{tr},*}(\scr F_n)$ define an object in $\mathsf{\Delta}{\scr Q}_{\mathrm{tr}}$.  
\end{lemma}

\noindent Now we have the following obvious lemma: 

\begin{lemma}\label{QSymbbZtr}
The Nisnevich sheaves with transfers $\bbZtr[\bbG_m^{\stimes n}]$ form a family of Nisnevich sheaves with transfers with a commutative external product. 
\end{lemma}

\begin{proof} Let $X,Y$ be schemes. The tensor product of correspondences gives us morphisms
\begin{equation*}
\xymatrix{{\Corrf(X,\bbG_m^{\times n})\otimes\Corrf(Y,\bbG_m^{\times m})}\ar[r]\ar@{=}[d] & {\Corrf(X\times Y,\bbG^{\times n+m}_m)}\ar@{=}[d]\\
{\bbZtr[\bbG_m^{\times n}](X)\otimes\bbZtr[\bbG_m^{\times m}](Y)}\ar[r] & {\bbZtr[\bbG^{\times n+m}](X\times Y)}}
\end{equation*}
which gives a commutative external product on the $\bbZtr[\bbG^{\times n}_m]$. It is enough to remark that we have the following commutative diagram 
\begin{equation*}
\xymatrix{{\begin{tabular}{c}${\cal D}_n(X)$\\$\otimes$\\${\cal D}_m(Y)$\end{tabular}}\ar[r]\ar[d] & {\begin{tabular}{c}$\bbZtr[\bbG^{\times n}_m](X)$\\$\otimes$\\$\bbZtr[\bbG^{\times n}_m](Y)$ \end{tabular}}\ar[r]\ar[d] & {\begin{tabular}{c}$\bbZtr[\bbG^{\stimes n}_m](X)$\\ $\otimes$\\ $\bbZtr[\bbG^{\stimes n}_m](Y)$\end{tabular}}\ar[r]\ar@{.>}[d] & {0}\\
{{\cal D}_{n+m}(X\times Y)}\ar[r] & {\bbZtr[\bbG_m^{\times n+m}](X\times Y)}\ar[r] & {\bbZtr[\bbG_m^{\stimes n+m}](X\times Y)}\ar[r] & {0}}
\end{equation*}
to see that the quotients $\bbZtr[\bbG^{\stimes n}_m]$ inherit a commutative external product. 
\end{proof}

\begin{proposition}
The functor $\Upsilon^{\mathrm{eff}}_{\mathrm{cs}}$ is a commutative external product. In other words for $\Gamma=X(n)$ and $\Delta=Y(m)$ there exists a canonical morphism
   \begin{equation}\label{MorPropQM}
	\Upsilon^{\mathrm{eff}}_{\mathrm{cs}}(\Gamma)\otimestr\Upsilon^{\mathrm{eff}}_{\mathrm{cs}}(\Delta)\xrightarrow{\boxtimes^{\Upsilon^{\mathrm{eff}}_{\mathrm{cs}}}_{\Gamma,\Delta}}\Upsilon^{\mathrm{eff}}_{\mathrm{cs}}(\Gamma\times\Delta)
\end{equation}
	 and those morphisms are associative, commutative and functorial.

\end{proposition}

\begin{proof}
Fix two smooth quasi-projective schemes $X,Y$ and two non negative integers $n,m$. It is enough to construct a natural map in $\mathsf{\Delta}^{\op}\mathsf{\Delta}\Ptr$
$$\Upsilon^{\mathrm{eff}}_{\mathrm{cs}}(X(n))\otimesptr\Upsilon^{\mathrm{eff}}_{\mathrm{cs}}(Y(m))\ra\Upsilon^{\mathrm{eff}}_{\mathrm{cs}}((X\times Y)(n+m))$$
and so  to construct, for all smooth quasi-projective scheme $Z$, a natural map in the category $\mathsf{\Delta}^{\op}\mathsf{\Delta}\Mod(\bbZ)$
\begin{equation}\label{MorIntPropQM}
\Big[\Upsilon^{\mathrm{eff}}_{\mathrm{cs}}(X(n))\otimesptr\Upsilon^{\mathrm{eff}}_{\mathrm{cs}}(Y(m))\Big](Z)\ra\Upsilon^{\mathrm{eff}}_{\mathrm{cs}}((X\times Y)(n+m))(Z)
\end{equation}
functorial for finite correspondences. By lemmas \ref{VarTensGodtr} and \ref{QSymbbZtr} we know that the family of cosimplicial Godement resolutions 
\begin{equation*}
{\scr H}^*_n=\scr G^{\tr,*}_{\Nis}(\bbZtr[\bbG^{\stimes n}_m])\qquad n\in\bbN
\end{equation*}
defines an object in $\mathsf{\Delta}\scr Q^{\tr}$. Therefore for all schemes $\scr X$ and $\scr Y$ we have a natural associative and commutative morphism that is functorial for finite correspondences 
\begin{equation}\label{MorFn}
{\scr H}^*_n(\scr X)\otimes {\scr H}^*_m(\scr Y)\ra {\scr H}^*_{n+m}(\scr X\times\scr Y).
\end{equation}
In particular we have associative and commutative morphisms in $\mathsf{\Delta}^{\op}\mathsf{\Delta}\Mod(\bbZ)$ that are functorial for finite correspondences
\begin{equation*}
\xymatrix{{{\scr H}^*_n(\Delta^*_{\scr X})\otimes {\scr H}^*_m(\Delta^*_{\scr Y})}\ar[r]^{\textrm{\raisebox{2.5ex}{\textrm{induced by (\ref{MorFn})}}}}\ar`d[rd][rd]_{\boxtimes^{{\scr H}^*}_{\scr X(n),\scr Y(m)}} & {{\scr H}^*_{n+m}(\Delta^*_{\scr X}\times\Delta^*_{\scr Y})}\ar[d]^{{\renewcommand{\arraystretch}{.5}\begin{tabular}{l}$\scriptstyle{\textrm{induced by the diagonal}}$\\$\scriptstyle{\textrm{immersion}}$\end{tabular}}}\\
{} & {{\scr H}^*_{n+m}(\Delta^*_{\scr X\times\scr Y}).}}
\end{equation*}
Besides if $\scr X$ and $\scr W$ are smooth quasi-projective schemes we have
\begin{equation*}
\begin{split}
\Upsilon^{\mathrm{eff}}_{\mathrm{cs}}(\scr X(n))(\scr W)&=\Hom(\bbZtr[\scr W],\Homotr(\bbZtr[\Delta^*_{\scr X}],{\scr H}^*_n))\\
&=\Hom(\bbZtr[\Delta^*_{\scr X\times\scr W}],{\scr H}^*_n)={\scr H}^*_n(\Delta^*_{\scr X\times\scr W}).
\end{split}
\end{equation*}
So we have the following commutative diagram
{\renewcommand{\arraystretch}{1.2}
\begin{eqnarray}\label{ProofPropQM}
\xymatrix@C=.22cm{{\bigoplus\limits_{W,W'\in\Sm_F}\begin{tabular}{l}$\Upsilon^{\mathrm{eff}}_{\mathrm{cs}}(X(n))(W)$\\ $\otimes\,\Upsilon^{\mathrm{eff}}_{\mathrm{cs}}(Y(m))(W')$ \\$\otimes\,\Corrf(Z,W\times W')$ \end{tabular}}\ar[r]^(.6){\textrm{\raisebox{7ex}{$\boxtimes^{{\scr H}^*}_{(X\times W)(n),(Y\times W')(m)}$}}}\ar@{->>}[d]^{(\ref{ProofPropQM})} & {\begin{tabular}{l}$\Upsilon^{\mathrm{eff}}_{\mathrm{cs}}((X\times Y)(n+m))(W\times W')$\\$\otimes\,\Corrf(Z,W\times W')$\end{tabular}}\ar[d]^{{{\renewcommand{\arraystretch}{.5}\begin{tabular}{l}$\scriptstyle{\textrm{structure of presheaf}}$\\$\scriptstyle{\textrm{with transfers}}$\end{tabular}}}}\\
{\Big[\Upsilon^{\mathrm{eff}}_{\mathrm{cs}}(X(n))\otimesptr\Upsilon^{\mathrm{eff}}_{\mathrm{cs}}(Y(m))\Big](Z)}\ar@{.>}[r] & {\Upsilon^{\mathrm{eff}}_{\mathrm{cs}}((X\times Y)(n+m))(Z)}}
\nonumber
\end{eqnarray}\addtocounter{equation}{1}}\noindent
in which the surjection (\ref{ProofPropQM}) is obtained from the formula of appendix \ref{AppendixNtr}. This gives the morphisms (\ref{MorIntPropQM}) we were looking for. It is a consequence of the construction that these morphisms are functorial for finite correspondences and that the induced morphisms (\ref{MorPropQM}) are associative and commutative.
\end{proof}

\noindent We may apply the construction given in paragraph \ref{ParPrincipeCons} to the commutative external product
\begin{equation*}
\xymatrix{{{\cal L}(F)^{\op}\times\bbN}\ar[r] & {\Sm_F^{\op}\times\bbN}\ar[r]^{\Upsilon^{\mathrm{eff}}_{\mathrm{cs}}} & {\mathsf{\Delta}^{\op}\mathsf{\Delta}\Ntr}}
\end{equation*}
to obtain a DG tensor functor
\begin{equation}\label{Upsilon2sh}
\Upsilon^{\mathrm{eff}\fksh}_2:{\cal A}_2^{\mathrm{eff}\fksh}(F)\ra\bC(\Ntr).
\end{equation}
and the following commutative diagram
\begin{equation*}
\xymatrix{{} & {} & {\bC(\Ntr)}\\
{} & {} & {\bCp(\Cm(\Ntr))}\ar[u]_{\mathrm{Tot}}\\
{{\cal A}^{\mathrm{eff}\fksh}_2(F)}\ar`u[uu][rruu]^{\Upsilon^{\mathrm{eff}\fksh}_2}\ar@/^2em/[rru]^{\left\lbrack\mathrm{cc}\Upsilon^{\mathrm{eff}}_1\right\rbrack^{\otimes,\fksh}}\ar[r]^{\mathfrak{c}} & {{\cal A}^{\mathrm{eff}}_2(F)} & {\Cp(\Cm(\Ntr))}\ar[u]_{\textrm{inclusion}}\\
{} & {{\cal A}^{\mathrm{eff}}_1(F)}\ar[u]^{\iota^c}\ar[r]^{\Upsilon^{\mathrm{eff}}_1}\ar[ru]^{\mathrm{cc}\Upsilon^{\mathrm{eff}}_1}\ar[lu]^{\iota^{\fksh}} & {\mathsf{\Delta}\Cm(\Ntr)}\ar[u]_{\mathrm{cc}}}
\end{equation*}
where  $\Upsilon^{\mathrm{eff}}_1$ denotes the canonical extension to ${\cal A}^{\mathrm{eff}}_1(F)$  of the functor
\begin{eqnarray*}
{\cal L}(F)^{\op}\times\bbN & \ra & \mathsf{\Delta}\Cm(\Ntr) \\
X(n)_f & \mapsto & \mathrm{c}\Upsilon^{\mathrm{eff}}_{\mathrm{cs}}(X(n))[-n].
\end{eqnarray*}
By construction the square
\begin{equation}\label{LemRep}
\begin{tabular}{c}
\begin{math}
\xymatrix{{\mathcal{A}_2^{\mathrm{eff}\fksh}(F)}\ar[r]^{\Upsilon_2^{\mathrm{eff}\fksh}} & {\bC(\Ntr)}\\
{\mathcal{L}(F)^{\op}\times\bbN}\ar[r]\ar[u] & {\Sm_F^{\op}\times\bbN}\ar[u]^{\Upsilon^{\mathrm{eff}}}}
\end{math}
\end{tabular}
\end{equation}
is commutative.

\subsubsection{Variant with support} 

We shall need during the construction a variant of the previous functor that takes support into account in order to use the semi-purity of motivic cohomology that assure the vanishing of some cohomology classes. For each objet $(X,f)\in{\cal L}(F)$ and each non negative integer $n$ we denote by $(X,f)^{(n)}$ the set of closed subschemes $W$ of $X$ such that
$$f^*(W):=X'\times_XW$$
is a pure codimension $n$ closed subscheme of $X'$.  Since $f$ has a smooth section $s$ the above condition implies that the closed subscheme 
$$W=s^*(f^*(W))$$
is also of pure codimension $n$ in $X$. The set $(X,f)^{(n)}$ is filtrant for the relation of inclusion and functorial. Indeed if $W$ belongs to $(X,f)^{(n)}$ and $p:(Y,g)\ra (X,f)$ is a morphism in ${\cal L}(F)$ then $p^*(W)$ belongs to  $(Y,g)^{(n)}$. This assertion is a consequence of the equality
$$g^*p^*(Z)=q^*f^*(Z)$$
where $q$ is a flat morphism that fits in a commutative square
\begin{equation*}
\xymatrix{{Y'}\ar[d]^{g}\ar[r]^{q} & {X'}\ar[d]^{f}\\
{Y}\ar[r]^{p} & {X.}}
\end{equation*}
In the sequel we denote by $\mathsf{R}\mathsf{\Gamma}_{\star}(X,\bbZ(n))$ the complex
\begin{equation*}
\mathsf{R}\mathsf{\Gamma}_{\star}(X,\bbZ(n)):=\colim_{W\in(X,f)^{(n)}}\mathsf{R}\mathsf{\Gamma}_{W}(X,\bbZ(n))
\end{equation*}
and its cohomology by $H^m_{\star}(X,\bbZ(n))_{\mathsf{V}}=H^m\mathsf{R}\mathsf{\Gamma}_{\star}(X,\bbZ(n))$. To define the variant with support we consider the functor
\begin{equation}\label{FoncSupp}
\Upsilon_{\mathrm{cs}}^{\mathrm{eff}\star}:{\mathcal{L}}(F)^{\op}\times\bbN\ra\mathsf{\Delta}^{\op}\mathsf{\Delta}\Ntr
\end{equation}
given for an object $(X,f)\in\mathcal{L}(F)$ and a non negative integer $n$ by 
\begin{equation*}
\Upsilon_{\mathrm{cs}}^{\mathrm{eff}\star}(X(n)_f)=\colim_{W\in (X,f)^{(n)}}\Upsilon^{\mathrm{eff}}_{W\mathrm{cs}}(X(n)).
\end{equation*}
  
\begin{proposition}
    The functor $\Upsilon_{\mathrm{cs}}^{\mathrm{eff}\star}$ is a commutative external product. In other words for $\Gamma=X(n)_f$ and $\Delta=Y(m)_g$ there exists a canonical morphism
    \begin{equation}\label{MorPropQMstar}
    \Upsilon_{\mathrm{cs}}^{\mathrm{eff}\star}(\Gamma)\otimestr\Upsilon_{\mathrm{cs}}^{\mathrm{eff}\star}(\Delta)\xrightarrow{\boxtimes^{\Upsilon_{\mathrm{cs}}^{\mathrm{eff}
    \star}}_{\Gamma,\Delta}}\Upsilon_{\mathrm{cs}}^{\mathrm{eff}\star}(\Gamma\times\Delta)
	\end{equation}
	 and those morphisms are associative, commutative and functorial. Futhermore the following square is commutative
	\begin{equation*}
	\xymatrix@C=2cm{{\Upsilon_{\mathrm{cs}}^{\mathrm{eff}\star}(\Gamma)\otimestr\Upsilon_{\mathrm{cs}}^{\mathrm{eff}\star}(\Delta)}\ar[r]^{\boxtimes^{\Upsilon_{\mathrm{cs}}^
    {\mathrm{eff}\star}}_{\Gamma,\Delta}}
	\ar[d] & {\Upsilon^{\mathrm{\mathrm{eff}\star}}_{\mathrm{cs}}(\Gamma\times\Delta)}\ar[d]\\
    {\Upsilon^{\mathrm{eff}}_{\mathrm{cs}}(\Gamma)\otimestr\Upsilon^{\mathrm{eff}}_{\mathrm{cs}}(\Delta)}\ar[r]^{\boxtimes^{\Upsilon^{\mathrm{eff}}_
    {\mathrm{cs}}}_{\Gamma,\Delta}} & 
	{\Upsilon^{\mathrm{eff}}_{\mathrm{cs}}(\Gamma\times\Delta).}}
	\end{equation*}
\end{proposition}

\begin{proof}
Given closed subschemes $W\in (X,f)^{(n)}$ and $Z\in (Y,g)^{(n)}$, we have a commutative diagram
\begin{equation*}
\xymatrix@C=2cm{{} & {0}\ar[d]\\
{\Upsilon^{\mathrm{eff}}_{\mathrm{cs}}(X(n))_W\otimestr\Upsilon^{\mathrm{eff}}_{\mathrm{cs}}(Y(m))_Z}\ar@{.>}[r]\ar[d] & {\Upsilon^{\mathrm{eff}}_{\mathrm{cs}}((X\times Y)(n+m))_{W\times Z}}\ar[d]\\
{\Upsilon^{\mathrm{eff}}_{\mathrm{cs}}(X(n))\otimestr\Upsilon^{\mathrm{eff}}_{\mathrm{cs}}(Y(m))}\ar[r]^{\boxtimes^{\Upsilon^{\mathrm{eff}}_{\mathrm{cs}}}_{X(n),Y(m)}}\ar[d] & {\Upsilon^{\mathrm{eff}}_{\mathrm{cs}}((X\times Y)(n+m))}\ar[d]\\
{\Upsilon^{\mathrm{eff}}_{\mathrm{cs}}(U(n))\otimestr\Upsilon^{\mathrm{eff}}_{\mathrm{cs}}(V(m))}\ar[r]^{\boxtimes^{\Upsilon^{\mathrm{eff}}_{\mathrm{cs}}}_{U(n),V(m)}} & {\Upsilon^{\mathrm{eff}}_{\mathrm{cs}}((U\times V)(n+m))}}
\end{equation*}
where $U$ is the open complement of $W$ and $V$ the open complement of $Z$. It is enough to remark that $Z\times W$ belongs to $(X\times Y,f\times g)^{(n+m)}$ and to take the colimit to obtain the morphism (\ref{MorPropQMstar}).
\end{proof}

\noindent Once again we may apply the construction given in paragraph \ref{ParPrincipeCons} to the commutative external product (\ref{FoncSupp}) to obtain a DG tensor functor
\begin{equation}\label{Upsilon2shstar}
\Upsilon^{\mathrm{eff}\fksh\star}_2:{\cal A}_2^{\mathrm{eff}\fksh}(F)\ra\bC(\Ntr)
\end{equation}
with a natural transformation 
\begin{equation*}
\xymatrix@C=2cm{{\cal A}_2^{\mathrm{eff}\fksh}(F)\rtwocell^{\raisebox{1em}{$\scriptstyle{\Upsilon^{\mathrm{eff}\fksh\star}_{2}}$}}_{\raisebox{-1em}{$\scriptstyle{\Upsilon^{\mathrm{eff}\fksh}_{2}}$}} & \bC(\Ntr)\\}.
\end{equation*} 
We have the following commutative diagram
\begin{equation*}
\xymatrix{{} & {} & {\bC(\Ntr)}\\
{} & {} & {\bCp(\Cm(\Ntr))}\ar[u]_{\mathrm{Tot}}\\
{{\cal A}^{\mathrm{eff}\fksh}_2(F)}\ar`u[uu][rruu]^{\Upsilon^{\mathrm{eff}\fksh\star}_2}\ar@/^2em/[rru]^{\left\lbrack\mathrm{cc}\Upsilon^{\mathrm{eff\star}}_1\right\rbrack^{\otimes,\fksh}}\ar[r]^{\mathfrak{c}} & {{\cal A}^{\mathrm{eff}}_2(F)} & {\Cp(\Cm(\Ntr))}\ar[u]_{\textrm{inclusion}}\\
{} & {{\cal A}_1(\Sm_k)}\ar[u]^{\iota^c}\ar[r]^{\Upsilon^{\mathrm{eff}\star}_1}\ar[ru]^{\mathrm{cc}\Upsilon^{\mathrm{eff}\star}_1}\ar[lu]^{\iota^{\fksh}} & {\mathsf{\Delta}\Cm(\Ntr)}\ar[u]_{\mathrm{cc}}}
\end{equation*}
where  $\Upsilon^{\mathrm{eff}\star}_1$ denotes the canonical extension to ${\cal A}^{\mathrm{eff}}_1(F)$  of the functor
\begin{eqnarray*}
{\cal L}(F)^{\op}\times\bbN & \ra & \mathsf{\Delta}\Cm(\Ntr) \\
X(n)_f & \mapsto & \mathrm{c}\Upsilon^{\mathrm{eff}\star}_{\mathrm{cs}}(X(n)_f)[-n].
\end{eqnarray*}
For $\Gamma=\bbZ_X(n)_f$ the isomorphisms (\ref{IsoCoh}) give an isomorphism in $\D(\bbZ)$
\begin{equation}\label{IsoCohSupp}
\mathsf{\Gamma}_{\mathrm{Nis}}(\Spec(F),\Upsilon_2^{\mathrm{eff}\fksh}(\Gamma))=\mathsf{R}\mathsf{\Gamma}_{\star}(X,\bbZ(n)).
\end{equation} 

\subsection{Cycle class maps}

Recall that in \cite{EqBlochVoe} Voevodsky has constructed for a smooth quasi-projective scheme $X$, a closed subscheme $W$ of $X$  and integers $n,m$ a natural isomorphism
\begin{equation}\label{VoeCycCl}
\CH^n_W(X,m)\xrightarrow{\fk{cl}_X^{n,m}}H^{2n-m}_W(X,\bbZ(n))_{\mathsf{V}}.
\end{equation}
In addition those morphisms are functorial, map intersection products to cup products and are also compatible with relativization long exact sequences.  

By the universal property of $\bbE$ \cite[Part II, Ch. II, 3.1.13]{MixMotLev} there exists a unique DG tensor functor
\begin{eqnarray*}
\bbE & \ra & \bC(\Ntr)\\
\mathfrak{e}^{\otimes n} & \mapsto & \bbZtr[F]^{\otimestr n}.
\end{eqnarray*}
Using this functor in order to extend the functors (\ref{Upsilon2sh}) and (\ref{Upsilon2shstar}) to the coproduct of DG categories ${\cal A}_2^{\mathrm{eff}\fksh}(F)[\bbE]$ we obtain a natural transformation of  DG tensor functors
\begin{equation}\label{FoncA2}
\xymatrix@C=2cm{{\cal A}_2^{\mathrm{eff}\fksh}(F)[\bbE]\rtwocell^{\raisebox{1em}{$\scriptstyle{\Upsilon^{\mathrm{eff}\fksh\star}_{2,\bbE}}$}}_{\raisebox{-1em}{$\scriptstyle{\Upsilon^{\mathrm{eff}\fksh}_{2,\bbE}}$}} & {\bC(\Ntr).}}
\end{equation}

\subsubsection{Adding cycle class maps}

In the sequel one sets for short $\bbZtr[F]=\bbZtr[\Spec(F)]$. From (\ref{IsoCohSupp}) we have a canonical surjection
\begin{equation}\label{SurjSupp}
Z^{2n}\bHom\Big(\bbZtr[F],\Upsilon^{\mathrm{eff}\fksh\star}_{2,\bbE}(\bbZ_X(n)_f)\Big)\twoheadrightarrow H^{2n}_{\star}(X,\bbZ(n))
\end{equation}
where $\bHom$ is the complex of morphism in the DG category $\bC(\Ntr)$. Thus for all non zero cycle $\alpha \in{\cal Z}^n(X)_f$, we may choose a morphism of degree $2n$
\begin{equation*}
\Upsilon^{\mathrm{eff}\fksh\star}_3([\alpha]):\bbZtr[F]\ra\Upsilon^{\mathrm{eff}\fksh\star}_{2,\bbE}(\bbZ_X(n)_f)
\end{equation*}
which maps to the cycle class $\fk{cl}_X^n([\alpha])$ through the surjection (\ref{SurjSupp}) and denote by $\Upsilon^{\mathrm{eff}\fksh}_3([\alpha])$ the morphism obtained after composition with the natural transformation
\begin{equation*}
\xymatrix@C=1.5cm{{\bbZtr[F]}\ar[r]_{\raisebox{-1.3em}{$\scriptstyle{\Upsilon^{\mathrm{eff}\fksh\star}_3([\alpha])}$}}\ar@/^2em/[rr]^{\Upsilon^{\mathrm{eff}\fksh}_3([\alpha])} & {\Upsilon^{\mathrm{eff}\fksh\star}_{2,\bbE}(\bbZ_X(n)_f)}\ar[r] & {\Upsilon^{\mathrm{eff}\fksh}_{2,\bbE}(\bbZ_X(n)_f).}}
\end{equation*}
Those choices define an extension of the DG tensor functors and the natural transformation (\ref{FoncA2}) to the DG tensor category without unit ${\cal A}_3^{\mathrm{eff}\fksh}(F)$
\begin{equation}\label{FoncA3}
\xymatrix@C=2cm{{\cal A}_3^{\mathrm{eff}\fksh}(F)\rtwocell^{\raisebox{1em}{$\scriptstyle{\Upsilon^{\mathrm{eff}\fksh\star}_3}$}}_{\raisebox{-1em}{$\scriptstyle{\Upsilon^{\mathrm{eff}\fksh}_3}$}} & {\bC(\Ntr).}}
\end{equation}
(a) For all morphism $p:(Y,g)\ra (X,f)$ in ${\cal L}(F)$ and all non zero cycle $\alpha\in{\cal Z}^n(X)_f$ the image of the closed morphism of degree $2n$
\begin{equation}\label{MorA4a}
p^*\circ[\alpha]-[p^*\alpha]
\end{equation}
by $\Upsilon^{\mathrm{eff}\fksh\star}_3$ is exact since it maps to 
\begin{equation*}
p^*\circ\fk{cl}_X^q([\alpha])-\fk{cl}_Y^q([p^*\alpha])=0
\end{equation*}  
in $H^{2n}_{\star}(X,\bbZ(n))$. So we may choose a morphism of degree $2n-1$
\begin{equation*}
\Upsilon^{\mathrm{eff}\fksh\star}_4(h_{X,Y,[\alpha],p^*}):\bbZtr[F]\ra \Upsilon^{\mathrm{eff}\fksh\star}_3(\bbZ_Y(n)_g)
\end{equation*}
whose differential is equal to the image by $\Upsilon^{\mathrm{eff}\fksh\star}_3$ of the morphism (\ref{MorA4a}). We denote $\Upsilon^{\mathrm{eff}\fksh}_4(h_{X,Y,[\alpha],p^*})$ the morphism obtained after composition with the natural transformation
\begin{equation*}
\xymatrix@C=2cm{{\bbZtr[F]}\ar[r]_{\raisebox{-1.3em}{$\scriptstyle{\Upsilon^{\mathrm{eff}\fksh\star}_4(h_{X,Y,[\alpha],p^*})}$}}\ar@/^2em/[rr]^{\Upsilon^{\mathrm{eff}\fksh}_4(h_{X,Y,[\alpha],p^*})} & {\Upsilon^{\mathrm{eff}\fksh\star}_{3}(\bbZ_X(n)_f)}\ar[r] & {\Upsilon^{\mathrm{eff}\fksh}_{3}(\bbZ_X(n)_f).}}
\end{equation*}
\par\noindent 
(b) For all non zero cycles $\alpha\in{\cal Z}^n(X)_f$ and $\beta\in {\cal Z}^m(Y)_g$, the images by $\Upsilon^{\mathrm{eff}\fksh\star}_3$ of the closed morphisms of degree $2(n+m)$
\begin{eqnarray}\label{MorA4b}
\boxtimes_{\Gamma,\Delta}\circ([\alpha]\otimes[\beta])-\boxtimes_{\Gamma\times\Delta,\bb1}\circ([\alpha\times\beta]\otimes[1])\\\label{MorA4c}
\boxtimes_{\Gamma,\Delta}\circ([\alpha]\otimes[\beta])-\boxtimes_{\bb1,\Gamma\times\Delta}\circ([1]\otimes[\alpha\times\beta])
\end{eqnarray}
where $\Gamma=\bbZ_X(n)_f$ and $\Delta=\bbZ_Y(m)_g$ are exact since they map to 
\begin{equation*}
\fk{cl}_X^n(\alpha)\otimes\fk{cl}_Y^{m}(\beta)-\fk{cl}_{X\times Y}^{n+m}(\alpha\times\beta)=0.
\end{equation*} So we may choose a morphism of degree  $2(n+m)-1$ 
\begin{eqnarray*}
\Upsilon^{\mathrm{eff}\fksh\star}_4(h^{l}_{X,Y,[\alpha],[\beta]}):\bbZtr[F]\ra \Upsilon^{\mathrm{eff}\fksh\star}_3(\Gamma\times\Delta)\\
\Upsilon^{\mathrm{eff}\fksh\star}_4(h^{r}_{X,Y,[\alpha],[\beta]}):\bbZtr[F]\ra \Upsilon^{\mathrm{eff}\fksh\star}_3(\Gamma\times\Delta)
\end{eqnarray*}
whose differentials are equal to the images by $\Upsilon^{\mathrm{eff}\fksh\star}_3$ of (\ref{MorA4b}) and (\ref{MorA4c}). We denote by $\Upsilon^{\mathrm{eff}\fksh}_4(h^{l}_{X,Y,[\alpha],[\beta]})$ and $\Upsilon^{\mathrm{eff}\fksh}_4(h^{r}_{X,Y,[\alpha],[\beta]})$ the morphisms obtained after composition with the natural transformation
\begin{equation*}
\xymatrix@C=2cm{{\bbZtr[F]}\ar[r]_{\raisebox{-1.3em}{$\scriptstyle{\Upsilon^{\mathrm{eff}\fksh\star}_4(h^{l}_{X,Y,[\alpha],[\beta]})}$}}\ar@/^2em/[rr]^{\Upsilon^{\mathrm{eff}\fksh}_4(h^{l}_{X,Y,[\alpha],[\beta]})} & {\Upsilon^{\mathrm{eff}\fksh\star}_{3}(\Gamma\times\Delta)}\ar[r] & {\Upsilon^{\mathrm{eff}\fksh}_{3}(\Gamma\times\Delta)}}
\end{equation*}
\begin{equation*}
\xymatrix@C=2cm{{\bbZtr[F]}\ar[r]_{\raisebox{-1.3em}{$\scriptstyle{\Upsilon^{\mathrm{eff}\fksh\star}_4(h^{r}_{X,Y,[\alpha],[\beta]})}$}}\ar@/^2em/[rr]^{\Upsilon^{\mathrm{eff}\fksh}_4(h^{r}_{X,Y,[\alpha],[\beta]})} & {\Upsilon^{\mathrm{eff}\fksh\star}_{3}(\Gamma\times\Delta)}\ar[r] & {\Upsilon^{\mathrm{eff}\fksh}_{3}(\Gamma\times\Delta).}}
\end{equation*}
\par\noindent 
(c) For all cycles $\alpha,\beta\in{\cal Z}^n(X)_f$ and all integers $a,b$ the image of the closed morphism of degree $2n$
\begin{equation}\label{MorA4d}
[a\alpha+b\beta]-a[\alpha]-b[\beta]
\end{equation}
by $\Upsilon^{\mathrm{eff}\fksh\star}_3$ is exact since it maps to 
\begin{equation*}
\fk{cl}_X^n(a\alpha+b\beta)-a\fk{cl}_X^{n}(\alpha)-b\fk{cl}_X^{n}(\beta)=0.
\end{equation*}  
in $H^{2n}_{\star}(X,\bbZ(n))$. So we may choose a morphism of degree  $2n-1$
\begin{equation*}
\Upsilon^{\mathrm{eff}\fksh\star}_4(h_{a,b,[\alpha],[\beta]}):\bbZtr[F]\ra \Upsilon^{\mathrm{eff}\fksh\star}_3(\bbZ_ {X}(n)_f)
\end{equation*}
whose differential is equal to the image by $\Upsilon^{\mathrm{eff}\fksh\star}_3$ of the morphism (\ref{MorA4d}). We denote by $\Upsilon^{\mathrm{eff}\fksh}_4(h_{a,b,[\alpha],[\beta]})$ the morphism obtained after composition with the natural transformation
\begin{equation*}
\xymatrix@C=2cm{{\bbZtr[F]}\ar[r]_{\raisebox{-1.3em}{$\scriptstyle{\Upsilon^{\mathrm{eff}\fksh\star}_4(h_{a,b,[\alpha],[\beta]})}$}}\ar@/^2em/[rr]^{\Upsilon^{\mathrm{eff}\fksh}_4(h_{a,b,[\alpha],[\beta]})} & {\Upsilon^{\mathrm{eff}\fksh\star}_{3}(\bbZ_X(n)_f)}\ar[r] & {\Upsilon^{\mathrm{eff}\fksh}_{3}(\bbZ_X(n)_f).}}
\end{equation*}
Those choices define an extension of the DG tensor functors and the natural transformation (\ref{FoncA3}) to the DG tensor category without unit ${\cal A}_4^{\mathrm{eff}\fksh}(F)$
\begin{equation}\label{FoncA4}
\xymatrix@C=2cm{{\cal A}_4^{\mathrm{eff}\fksh}(F)\rtwocell^{\raisebox{1.3em}{$\scriptstyle{\Upsilon^{\mathrm{eff}\fksh\star}_4}$}}_{\raisebox{-1.3em}{$\scriptstyle{\Upsilon^{\mathrm{eff}\fksh}_4}$}} & {\bC(\Ntr).}}
\end{equation}

\subsubsection{Cleaning up to homotopy}

Let $r\geq 1$ and $k\geq 1$ be integers. Suppose that the DG tensor functors and the natural transformation (\ref{FoncA4}) have been extended to the DG tensor category without unit ${\cal A}_4^{\mathrm{eff}\fksh}(F)^{(r,k-1)}$
\begin{equation}\label{FoncA4rka}
\xymatrix@C=2cm{{\cal A}_4^{\mathrm{eff}\fksh}(F)^{(r,k-1)}\rtwocell^{\raisebox{1.3em}{$\scriptstyle{\Upsilon^{\mathrm{eff}\fksh\star}_{4,(r,k-1)}}$}}_{\raisebox{-1.3em}{$\scriptstyle{\Upsilon^{\mathrm{eff}\fksh}_{4,(r,k-1)}}$}} & {\bC(\Ntr).}}
\end{equation}
Fix a closed morphism of degree $2n-r$ in ${\cal A}_4^{\mathrm{eff}\fksh}(F)^{(r,k-1)}$
\begin{equation}\label{MorNet}
s:\fk{e}^{\otimes k}\ra\bbZ_X(n)_f.
\end{equation}
Its image by the functor $\Upsilon^{\mathrm{eff}\fksh\star}_{4,(r,k-1)}$ is exact since by semi-purity of motivic cohomology we have
$$H^{2n-r}_{\star}(X,\bbZ(n))=0.$$
So we may choose a morphism of degree $2n-r-1$ 
\begin{equation*}
\Upsilon^{\fksh,(*)}_{4,(r,k)(h_s)}:\bbZtr[F]\ra \Upsilon^{\mathrm{eff}\fksh\star}_{4,(r,k-1)}(\bbZ_X(n)_f)
\end{equation*}
whose differential is equal to the image by $\Upsilon^{\mathrm{eff}\fksh\star}_{4,(r,k-1)}$ of the morphism $s$ and consider the morphism gotten after the composition with the natural transformation
\begin{equation*}
\xymatrix@C=1.5cm{{\bbZtr[F]}\ar[r]_{\raisebox{-1.3em}{$\scriptstyle{\Upsilon^{\mathrm{eff}\fksh\star}_{4,(r,k)}(h_s)}$}}\ar@/^2em/[rr]^{\Upsilon^{\mathrm{eff}\fksh}_{4,(r,k})(h_s)} & {\Upsilon^{\mathrm{eff}\fksh\star}_{4,(r,k-1)}(\bbZ_X(n)_f)}\ar[r] & {\Upsilon^{\mathrm{eff}\fksh}_{4,(r,k-1)}(\bbZ_X(n)_f).}}
\end{equation*}
\noindent Those choices define an extension of the DG tensor functors and the natural transformation (\ref{FoncA4rka}) to the DG tensor category without unit ${\cal A}_4^{\mathrm{eff}\fksh}(F)^{(r,k)}$
\begin{equation}\label{FoncA4rkb}
\xymatrix@C=2cm{{\cal A}_4^{\mathrm{eff}\fksh}(F)^{(r,k)}\rtwocell^{\raisebox{1.3em}{$\scriptstyle{\Upsilon^{\mathrm{eff}\fksh\star}_{4,(r,k)}}$}}_{\raisebox{-1.3em}{$\scriptstyle{\Upsilon^{\mathrm{eff}\fksh}_{4,(r,k)}}$}} & {\bC(\Ntr).}}
\end{equation}
By induction this extends the DG tensor functors and the natural transformation (\ref{FoncA4rkb}) to the DG tensor category without unit ${\cal A}_5^{\mathrm{eff}\fksh}(F)$
\begin{equation}\label{FoncA5}
\begin{tabular}{c}
\begin{math}
\xymatrix@C=2cm{{\cal A}_5^{\mathrm{eff}\fksh}(F)\rtwocell^{\raisebox{1.3em}{$\scriptstyle{\Upsilon^{\mathrm{eff}\fksh\star}_5}$}}_{\raisebox{-1.3em}{$\scriptstyle{\Upsilon^{\mathrm{eff}\fksh}_5}$}}\ar@{=}[d] & {\bC(\Ntr).}\\
{\twocolim_{r,k}{\cal A}_4^{\mathrm{eff}\fksh}(F)^{(r,k)}} & {}}
\end{math}
\end{tabular}
\end{equation}
We denote by $\Upsilon^{\mathrm{eff}\fksh}_{\textrm{mot}}$ the restriction of the functor $\Upsilon^{\mathrm{eff}\fksh}_5$  to the DG tensor category without unit ${\cal A}^{\mathrm{eff}}_{\mathrm{mot}}(F)$.

\subsection{The motivic commutative external product}

\noindent By applying the functor $\bCb$ and composing with the DG functor $\mathsf{Tot}$ one obtains a DG tensor functor
\begin{equation*}
\xymatrix@C=1.5cm{{\bC^{\mathrm{b},\mathrm{eff}}_{\mathrm{mot}\fksh}(F)}\ar[r]^{\bCb(\Upsilon^{\mathrm{eff}\fksh}_{\textrm{mot}})}\ar@/^3em/[rr]^{\Upsilon_{\mathrm{\bf{mot}}\fksh}^{\mathrm{eff}+}} & {\bCb(\bC(\Ntr))}\ar[r]^{\mathsf{Tot}} &  {\bC(\Ntr).}}
\end{equation*}
Passing to the homotopy categories gives a triangulated tensor functor
\begin{equation*}
\xymatrix@C=1.5cm{{\Kho^{\mathrm{b},\mathrm{eff}}_{\mathrm{mot}\fksh}(F)}\ar[r]_{\mathsf{Ho}(\Upsilon_{\mathrm{\bf{mot}}\fksh}^{\mathrm{eff}+})}\ar@/^2em/[rr]^{\Upsilon_{\mathrm{mot}\fksh}^{\mathrm{eff}}} & {\Kho(\Ntr)\ar[r]
} & {\D(\Ntr)}}
\end{equation*}
The square (\ref{LemRep}) and proposition \ref{lemGodRe} give a natural isomorphism
\begin{equation*}
\Upsilon^{\mathrm{eff}}_{\mathrm{mot}\fksh}(\Gamma)=\Homi(M(X),\bbZ(n))\qquad \Gamma=\bbZ_X(n)_f
\end{equation*}
which belongs to $\D^{-}(\Ntr)$. Since $\Kho^{\mathrm{b},\mathrm{eff}}_{\mathrm{mot}\fksh}(F)$ is generated as tensor triangulated category by those $\Gamma$ it follows that the functor $\Upsilon_{\mathrm{mot}\fksh}^{\mathrm{eff}}$ takes its values in $\D^{-}(\Ntr)$ and so induces a triangulated tensor functor
\begin{equation*}
\Upsilon_{\mathrm{mot}\fksh}^{\mathrm{eff}}:\Kho^{\mathrm{b},\mathrm{eff}}_{\mathrm{mot}\fksh}(F)\ra DM^{\mathrm{eff}}_{-}(F)
\end{equation*}
and by construction the square below is $2$-commutative
\begin{equation}\label{FoncRepDual}
\begin{tabular}{c}
\begin{math}
\xymatrix@C=1.5cm{{\Kho^{\mathrm{b},\mathrm{eff}}_{\mathrm{mot}\fksh}(F)}\ar[r]^{\Upsilon_{\mathrm{mot}\fksh}^{\mathrm{eff}}} & {DM^{\mathrm{eff}}_{-}(F)}\\
{{\cal L}(F)^{\op}\times\bbN}\ar[r]\ar[u]^{X(n)_f\mapsto\bbZ_X(n)_f} & {\Sm_F^{\op}\times\bbN.}\ar[u]_{X(n)\mapsto \Homi(M(X),\bbZ(n))}}
\end{math}
\end{tabular}
\end{equation}

\begin{proposition}\label{Upsilonsheff}
The functor $\Upsilon_{\mathrm{mot}\fksh}^{\mathrm{eff}}$ extends to ${\cal DM}^{\mathrm{eff}}_{\fksh}(F)$, in other words we have a $2$-commutative diagram
\begin{equation*}
\xymatrix{{\Kho^{\mathrm{b},\mathrm{eff}}_{\mathrm{mot}\fksh}(F)}\ar[r]\ar[rrd]_{\Upsilon_{\mathrm{mot}\fksh}^{\mathrm{eff}}} & {\D^{\mathrm{b},\mathrm{eff}}_{\mathrm{mot}\fksh}(F)}\ar[r]\ar@{.>}[rd] & {{\cal DM}^{\mathrm{eff}}_{\fksh}(F)}\ar@{.>}[d]^{\Upsilon_{\fksh}^{\mathrm{eff}}}\\
{} & {} & {DM^{\mathrm{eff}}_{-}}(F)}
\end{equation*}
where $\Upsilon_{\fksh}^{\mathrm{eff}}$ is a triangulated commutative external product.
\end{proposition}

\begin{proof}
Consider the commutative external product
\begin{equation*}
\Kho^{\mathrm{b},\mathrm{eff}}_{\mathrm{mot}\fksh}(F)^*\xrightarrow{i_{\mathrm{mot}}}\Kho^{\mathrm{b},\mathrm{eff}}_{\mathrm{mot}\fksh}(F)\xrightarrow{\Upsilon_{\mathrm{mot}\fksh}^{\mathrm{eff}}}DM^{\mathrm{eff}}_{-}(F)
\end{equation*}
which shall be denoted by $\Upsilon_{\mathrm{mot}\fksh}^{\mathrm{eff}*}$ in the sequel of this proof. It is enough to check that it induces a commutative external product
\begin{equation*}
\D^{\mathrm{b},\mathrm{eff}}_{\mathrm{mot}\fksh}(F)^*\xrightarrow{\Upsilon_{\mathrm{mot}\fksh}^{\mathrm{eff}*}} DM^{\mathrm{eff}}_{-}(F)
\end{equation*}
since in that case we have a $2$-commutative diagram
\begin{equation*}
\xymatrix{{\Kho^{\mathrm{b},\mathrm{eff}}_{\mathrm{mot}\fksh}(F)}\ar[r]\ar@/^2em/[rr]^{\Upsilon_{\mathrm{mot}\fksh}^{\mathrm{eff}}}\ar[d]^{r_{\mathrm{mot}}} & {\D^{\mathrm{b},\mathrm{eff}}_{\mathrm{mot}\fksh}(F)}\ar[d]^{r_{\mathrm{mot}}} & {DM^{\mathrm{eff}}_{-}(F)}\\
{\Kho^{\mathrm{b},\mathrm{eff}}_{\mathrm{mot}\fksh}(F)^*}\ar[r] & {\D^{\mathrm{b},\mathrm{eff}}_{\mathrm{mot}\fksh}(F)^*.}\ar[ru]_{\Upsilon_{\mathrm{mot}\fksh}^{\mathrm{eff}*}} & {}}
\end{equation*}
(a) {\itshape{Homotopy:}}
Let $p:(X,f)\ra (Y,g)$ be a morphism in ${\cal L}(F)$ such that $p:X\hookrightarrow Y$ is the inclusion of a pure codimension $1$ closed subscheme. Let $\hat{Y}$ be a closed subscheme of $Y$ and suppose that the closed subscheme $\hat{X}=X\times_Y\hat{Y}$ of $X$ is smooth and that we have a commutative diagram
\begin{equation*}
\xymatrix{{\hat{X}}\ar[r] & {\hat{Y}}\\
{\hat{X}\times\{0\}}\ar[r]\ar@{=}[u] & {\hat{X}\times\bbA^1.}\ar@{-}[u]_{\mathrm{iso.}}}
\end{equation*}
The square (\ref{FoncRepDual}) gives a commutative square
\begin{equation*}
\xymatrix@C=1.5cm{{\Upsilon_{\mathrm{mot}\fksh}^{\mathrm{eff}*}(\bbZ_{Y,\hat{Y}}(n)_g)}\ar[r]^{\Upsilon_{\mathrm{mot}\fksh}^{\mathrm{eff}*}(p^*)} & {\Upsilon_{\mathrm{mot}\fksh}^{\mathrm{eff}*}(\bbZ_{X,\hat{X}}(n)_f)}\\
{\Homi(M_{\hat{Y}}(Y),\bbZ(n))}\ar@{-}[u]^{\mathrm{iso.}}\ar[r]^{-\circ M(p)} & {\Homi(M_{\hat{X}}(X),\bbZ(n)).}\ar@{-}[u]^{\mathrm{iso.}}}
\end{equation*}
Using Gysin isomorphisms and $\bbA^1$-invariance one has $M_{\hat{Y}}(Y)=M_{\hat{X}}(X)$ and so the image of the morphism 
\begin{equation*}
p^*:\bbZ_{Y,\hat{Y}}(n)_g\ra\bbZ_{X,\hat{X}}(n)_{f}
\end{equation*}
by the functor $\Upsilon_{\mathrm{mot}\fksh}^{\mathrm{eff}*}$ is an isomorphism in $DM^{\mathrm{eff}}_{-}(F)$.
\par\medskip\noindent
(b) {\itshape{Excision:}} Let $(X,f)$ be an object of ${\cal L}(F)$, $\hat{X}$ a closed subscheme of $X$ and $U$ an open subset of $X$ containing $\hat{X}$. The square (\ref{FoncRepDual}) gives a commutative diagram 
\begin{equation*}
\xymatrix@C=1.5cm{{\Upsilon_{\mathrm{mot}\fksh}^{\mathrm{eff}*}(\bbZ_{X,\hat{X}}(n)_f)}\ar[r]^{\Upsilon_{\mathrm{mot}\fksh}^{\mathrm{eff}*}(j^*)}\ar@{-}[d]^{\mathrm{iso.}} & {\Upsilon_{\mathrm{mot}\fksh}^{\mathrm{eff}*}(\bbZ_{U,\hat{X}}(n)_{j^*f})}\ar@{-}[d]^{\mathrm{iso.}}\\
{\Homi(M_{\hat{X}}(X),\bbZ(n))}\ar[r]^{-\circ M(j)} & {\Homi(M_{\hat{X}}(U)\bbZ(n)).}}
\end{equation*}
By excision we have $M_{\hat{X}}(X)=M_{\hat{X}}(U)$ and so the image of the morphism 
\begin{equation*}
j^*:\bbZ_{X,\hat{X}}(n)_f\ra\bbZ_{U,\hat{X}}(n)_{j^*f}
\end{equation*}
by the functor $\Upsilon_{\mathrm{mot}\fksh}^{\mathrm{eff}*}$ is an isomorphism in $DM^{\mathrm{eff}}_{-}(F)$.
\par\medskip\noindent
(d) {\itshape{Gysin isomorphism:}} Let $p:(P,g)\ra (X,f) $ be a morphism in ${\cal L}(F)$ such that the morphism $p:P\ra X$ is a smooth morphism that has a section $s$ which is a pure codimension $d$ closed immersion. Assume that $[X]$ belongs to ${\cal Z}^d(P)_g$. It is enough to check that the image by the functor $\Upsilon_{\mathrm{mot}\fksh}^{\mathrm{eff}*}$ of the morphism 
\begin{equation*}
\xymatrix{{\fk{e}\otimes\bbZ_X(n-d)_f[-2d]\oplus\bbZ_{P\times P,X\times P}(n)_{g\times g\amalg\Delta}}\ar[d]_{{\left(\!\sideset{_0^\alpha}{_{\Delta^*}^{-\rho}}{\mathop{\,}}\right)}}\\{\bbZ_{P\times P,X\times P}(n)_{g\times g}\oplus\bbZ_{P,X}(n)_g}}
\end{equation*}
is an isomorphism. This amounts to check that the mapping cone of the morphism 
\begin{equation*}
\Upsilon_{\mathrm{\bf{mot}}\fksh}^{\mathrm{eff}*}\begin{pmatrix} \alpha & -\rho\\ 0 & \Delta^*\end{pmatrix}=\begin{pmatrix}\Upsilon_{\mathrm{\bf{mot}}\fksh}^{\mathrm{eff}*}(\alpha) & -\Id \\ 0 & \Upsilon_{\mathrm{\bf{mot}}\fksh}^{\mathrm{eff}*}(\Delta^*)\end{pmatrix}
\end{equation*}
in $\mathrm{C}(\Ntr)$ is isomorphic to zero in $DM^{\mathrm{eff}}_{-}(F)$. So it is enough to check\footnote{We use the following remark. Let $\mathsf{A}$ be an additive category, $A\xrightarrow{a}B\xrightarrow{b}C $ be morphisms in $\mathrm{C}(\mathsf{A})$. The canonical morphism between the mapping cones defined by the diagram 
\begin{equation*}
\xymatrix{{A}\ar[r]^{b\circ a}\ar[d]_{(\Id,a)} & {C}\ar[d]^{(0,\Id)}\ar[r] & {\mathsf{Mc}(b\circ a)}\ar@{.>}[d]\\
{A\oplus B}\ar[r]^{{\left(\!\sideset{_0^a}{_{\;\;b}^{-\Id}}{\mathop{\,}}\right)}} & {B\oplus C}\ar[r] & {\mathsf{Mc}\begin{pmatrix}a & -\Id\\ 0& b\end{pmatrix}}}
\end{equation*}
is a quasi-isomorphism. } that the mapping cone of the morphism 
$$\Upsilon_{\mathrm{\bf{mot}}\fksh}^{\mathrm{eff}*}(\Delta^*)\circ\Upsilon_{\mathrm{\bf{mot}}\fksh}^{\mathrm{eff}*}(\alpha) $$
is isomorphic to zero in $DM^{\mathrm{eff}}_{-}(F)$, which means that the morphism
$$\Upsilon_{\mathrm{mot}\fksh}^{\mathrm{eff}*}(\Delta^*)\circ\Upsilon_{\mathrm{mot}\fksh}^{\mathrm{eff}*}(\alpha) $$
is an isomorphism $DM^{\mathrm{eff}}_{-}(F)$. The square (\ref{FoncRepDual}) gives a commutative diagram 
\begin{equation}\label{MorCanc}
\xymatrix{{\Upsilon_{\mathrm{mot}\fksh}^{\mathrm{eff}*}(\fk{e}\otimes\bbZ_X(n-d)_f[-2d])}\ar@{-}[d]^{\mathrm{iso.}}\ar[r]^{\raisebox{1.3em}{$\scriptstyle{\Upsilon_{\mathrm{mot}\fksh}^{\mathrm{eff}*}(\Delta^*)\circ\Upsilon_{\mathrm{mot}\fksh}^{\mathrm{eff}*}(\alpha)}$}} & {\Upsilon_{\mathrm{mot}\fksh}^{\mathrm{eff}*}(\bbZ_{P,X}(n)_g)}\ar@{-}'`r[r]'`[rddd]^{\mathrm{iso.}}[ddd]  & *\txt<0pc>{}\\
{\Homi(M(X),\bbZ(n-d))[-2d]}\ar[r]^{\raisebox{2.3em}{$\scriptstyle\big[\fk{cl}^d_P([X])[-2d]\big]\otimes\big[-\circ M(p)\big]$}}\ar@{=}[d] & {\begin{tabular}{c}$\Homi(M_X(P),\bbZ(d))$\\$\otimes$\\$\Homi(M(P),\bbZ(n-d))$\end{tabular}}\ar[d]^{\boxtimes} & {}\\
{\Homi(M(X)[2d],\bbZ(n-d))}\ar[d]^{(\ref{MorCanc})} & {\Homi(M_{X\times P}(P\times P),\bbZ(n))}\ar[d]^{-\circ M(\Delta)} & *\txt<0pc>{}\\
{\Homi(M(X)(d)[2d],\bbZ(n))}\ar[r]_{\raisebox{-1.3em}{$\scriptstyle-\circ\big[\fk{cl}^d_{P}([X])\cupp M(p)\big]$}} & {\Homi(M_X(P),\bbZ(n))} & {}}
\nonumber
\end{equation}\addtocounter{equation}{1}\noindent
and the result follows from the fact that the Gysin isomorphism $s^*:M_X(P)\ra M(X)(d)[2d]$ is given by\footnote{For this classical result we refer for example to \cite[Remark 2.3]{IvoB}.}
\begin{equation*}
s^*=\fk{cl}^d_{P}([X])\cupp M(p)
\end{equation*}
and from Voevodsky's cancellation theorem \cite{Cancel} which shows that the map (\ref{MorCanc}) is an isomorphism.
\par\medskip\noindent
(e) {\itshape{Moving lemma:}} Let $(X,f)$ be an object in ${\cal L}(F)$ and $g:Z\ra X$ a morphism in $\Sm_F$.  To check that the image by $\Upsilon_{\mathrm{mot}\fksh}^{\mathrm{eff}*}$ of the morphism $\rho_{f,g}$ is an isomorphism, it is enough to remark that by construction we have a commutative square
\begin{equation*}
\xymatrix@C=2cm{{\Upsilon_{\mathrm{mot}\fksh}^{\mathrm{eff}*}(\bbZ_X(n)_f)}\ar[r]^{\Upsilon_{\mathrm{mot}\fksh}^{\mathrm{eff}*}(\rho_{f,g})} & {\Upsilon_{\mathrm{mot}\fksh}^{\mathrm{eff}*}(\bbZ_X(n)_{f\amalg g})}\\
{\mathsf{Tot}\Big[\mathrm{cc}\Big(\mathrm{c}\Upsilon^{\mathrm{eff}}_{\mathrm{cs}}(X(n))[-n]\Big)\Big]}\ar@{=}[r]\ar@{=}[u] & {\mathsf{Tot}\Big[\mathrm{cc}\Big(\mathrm{c}\Upsilon^{\mathrm{eff}}_{\mathrm{cs}}(X(n))[-n]\Big)\Big].}\ar@{=}[u]}
\end{equation*} 
\par\medskip\noindent
(f) {\itshape{Unit:}} We have in $DM^{\mathrm{eff}}_{-}(F)$ a commutative diagram
\begin{equation*}
\xymatrix@C=2cm{{\Upsilon_{\mathrm{mot}\fksh}^{\mathrm{eff}*}(\fk{e}\otimes\bb1)}\ar[r]^{\Upsilon_{\mathrm{mot}\fksh}^{\mathrm{eff}*}([[1]]\otimes\Id)} & {\Upsilon_{\mathrm{mot}\fksh}^{\mathrm{eff}*}(\bb1\otimes\bb1)}\\
{\Upsilon_{\mathrm{mot}\fksh}^{\mathrm{eff}*}(\fk{e})\otimes\Upsilon_{\mathrm{mot}\fksh}^{\mathrm{eff}*}(\bb1)}\ar[u]_{\mathrm{iso.}}\ar@{-}[d]_{\mathrm{iso.}}\ar[r]^{\Upsilon_{\mathrm{mot}\fksh}^{\mathrm{eff}*}([[1]])\otimes\Id} & {\Upsilon_{\mathrm{mot}\fksh}^{\mathrm{eff}*}(\bb1)\otimes\Upsilon_{\mathrm{mot}\fksh}^{\mathrm{eff}*}(\bb1)}\ar@{-}[d]_{\mathrm{iso.}}\ar[u]_{\mathrm{iso.}}\\
{\bbZ\otimes\bbZ}\ar[r]_{\fk{cl}^0_{\Spec(F)}([1])\otimes\Id} & {\bbZ\otimes\bbZ}}
\end{equation*}
and so it is enough to remark that $\fk{cl}^0_{\Spec(F)}([1])=\Id_{\bbZ}$ to see that the image by the functor $\Upsilon_{\mathrm{mot}\fksh}^{\mathrm{eff}*}$  of the morphism 
\begin{equation*}
[[1]]\otimes\Id:\fk{e}\otimes\bb1\ra\bb1\otimes\bb1
\end{equation*}
is an isomorphism in $DM^{\mathrm{eff}}_{-}(F)$.
\end{proof}

\begin{remark}\label{CondTens}
 Recall that $\Upsilon^{\mathrm{eff}}_{\mathrm{mot}\fksh}$ and $r_{\mathrm{mot}}$ are tensor functors. By construction the commutative external product on $\Upsilon^{\mathrm{eff}}_{\fksh}$  is given by the composition
\begin{equation*}
\xymatrix{{\Upsilon^{\mathrm{eff}}_{\fksh}(\Gamma)\otimes\Upsilon^{\mathrm{eff}}_{\fksh}(\Delta)}\ar[ddd]\ar@{=}[r] & {\Big[\Upsilon^{\mathrm{eff}}_{\mathrm{mot}\fksh}\circ i_{\mathrm{mot}}\circ r_{\mathrm{mot}}(\Gamma)\Big]\otimes\Big[\Upsilon^{\mathrm{eff}}_{\mathrm{mot}\fksh}\circ i_{\mathrm{mot}}\circ r_{\mathrm{mot}}(\Delta)\Big]}\ar[d]^{\mathrm{iso.}}\\
{} & {\Upsilon^{\mathrm{eff}}_{\mathrm{mot}\fksh}\left(\Big[i_{\mathrm{mot}}\circ r_{\mathrm{mot}}(\Gamma)\Big]\otimes\Big[ i_{\mathrm{mot}}\circ r_{\mathrm{mot}}(\Delta)\Big]\right)}\ar[d]^{\Upsilon^{\mathrm{eff}}_{\mathrm{mot}\fksh}(\boxtimes_{\mathrm{mot}})}\\
{} & {\Upsilon^{\mathrm{eff}}_{\mathrm{mot}\fksh}\circ i_{\mathrm{mot}}\Big(r_{\mathrm{mot}}(\Gamma)\times r_{\mathrm{mot}}(\Delta)\Big)}\ar[d]^{\mathrm{iso.}}\\
{\Upsilon^{\mathrm{eff}}_{\fksh}(\Gamma\otimes\Delta)}\ar@{=}[r] & {\Upsilon^{\mathrm{eff}}_{\mathrm{mot}\fksh}\circ i_{\mathrm{mot}}\circ r_{\mathrm{mot}}(\Gamma\otimes\Delta)}}
\end{equation*}
for two objects $\Gamma,\Delta$ in $\D^{\mathrm{b},\mathrm{eff}}_{\mathrm{mot},\fksh}(F)$. Since this category is generated as a triangulated category by the motives $\Gamma=\bbZ_X(n)$ the commutative external product $\Upsilon^{\mathrm{eff}}_{\fksh}$ is a tensor functor if the image by $\Upsilon^{\mathrm{eff}}_{\mathrm{mot}\fksh}$ of the K\"unneth map
\begin{equation*}
\boxtimes_{\mathrm{mot}}=\boxtimes^{\fksh}_{\Gamma,\Delta}:\Gamma\otimes\Delta\ra\Gamma\times\Delta
\end{equation*}
is an isomorphism where $\Gamma=\bbZ_X(n)$ and $\Delta=\bbZ_Y(m)$. 
\end{remark}

\subsection{Comparison of motivic cohomologies}

Let $\Gamma$ be an object in $\mathcal{DM}^{\mathrm{eff}}_{\fksh}(F)$ and $n$ a non negative integer. We are going to check that the induced map 
\begin{equation*}
\Hom_{\mathcal{DM}^{\mathrm{eff}}_{\fksh}(F)}\big(\bbZ(n),\Gamma\big)\xrightarrow{\Upsilon^{\mathrm{eff}}_{\fksh}}\Hom_{DM^{\mathrm{eff}}_{-}(F)}\big(\bbZ(n), \Upsilon^{\mathrm{eff}}_{\fksh}(\Gamma)\big)
\end{equation*}
is an isomorphism.
Let $X$ be a smooth quasi-projective scheme and $n$ a non negative integer. In particular we have a map
\begin{equation*}
\xymatrix@C=.2cm{{\Hom_{\mathcal{DM}^{\mathrm{eff}}_{\fksh}(F)}(\bb1,\bbZ_X(n)[m])}\ar@{=}[dd]\ar[r] & {\Hom_{DM^{\mathrm{eff}}_{-}(F)}(\bbZ,\Homi(M(X),\bbZ(n))[m])}\ar@{=}[d]\\
{} & {\Hom_{DM^{\mathrm{eff}}_{-}(F)}(M(X),\bbZ(n)[m])}\ar@{=}[d]\\
{H^m(X,\bbZ(n))_{\mathsf{L}}}\ar[r]^{\Upsilon_{\fksh}^{\mathrm{eff}}} & {H^m(X,\bbZ(n))_{\mathsf{V}}}}
\end{equation*}
on the motivic cohomologies. As both motivic cohomologies are isomorphic to higher Chow groups through respective cycle class maps, to check that this map is an isomorphism one has just to check its compatibility with cycle class maps. This is done in proposition \ref{CompMotCoh} and for this we give in lemma \ref{LemDescCyclLev} a basic property of Levine's cycle class map.\par
Recall that in \cite[Part I, Ch. I, 3.3.1.ii]{MixMotLev} the cycle functor (\ref{ElCycl}) is extended to a DG functor
\begin{equation}\label{Zmot}
\mathcal{Z}_{\mathrm{mot}}:\bC^{\mathrm{b}}_{\mathrm{mot}}(F)^*\ra\bC^{\mathrm{b}}(\bbZ)
\end{equation}
which satisfies the two following main properties:
\begin{itemize}
\item{The maps (a), (b), (c) of step 3 are sent to zero.}
\item{For an object $(X,f)\in\mathcal{L}(F)$ and a non-negative integer $n$ one has
\begin{equation*}
\mathcal{Z}_{\mathrm{mot}}(\bbZ_X(n)_f[2n])=\mathcal{Z}^{n}(X)_f.
\end{equation*}}
\end{itemize}
Since its codegeneracy maps are flat, the cosimplicial scheme $\Delta^*$ is a \textit{very smooth} cosimplicial scheme in the sense of \cite[Part I, Ch. I, \S2.4.1]{MixMotLev}. So it may be lifted to a cosimplicial object of ${\cal L}(F)$
\begin{equation*}
(\Delta^*,\delta^*):\mathsf{\Delta}^n\mapsto (\Delta^n,\delta^n)
\end{equation*}  
where
\begin{equation*}
{\Delta^n}':=\coprod_{\mathsf{\Delta}^n_{m,\mathrm{nd}}}\Delta^{m}
\end{equation*}
and $\delta^n:{\Delta^n}'\ra\Delta^n$ is the map given by $\Delta(s)$ on the component indexed by the non-degenerate $m$-simplex $s$ of $\mathsf{\Delta}^n$. For an object $\Gamma$ in $\bC^{\mathrm{b}}_{\mathrm{mot}}(F)^*$ this gives a cosimplicial object of $\mathsf{Z^0}\bC^{\mathrm{b}}_{\mathrm{mot}}(F)^*$
\begin{equation*}
\Gamma\times\bbZ_{\Delta^*}(0)_{\delta^*}:\mathsf{\Delta}^n\mapsto \Gamma\times\bbZ_{\Delta^n}(0)_{\delta^n}
\end{equation*}
and the Bloch's cycle complex attached to $\Gamma$ is given by  
\begin{equation*}
{\cal Z}_{\mathrm{mot}}(\Gamma,*):=\mathsf{Tot}\Big[\mathrm{cc}{\cal Z}_{\mathrm{mot}}\big(\Gamma\times\bbZ_{\Delta^*}(0)_{\delta^*}\big)\Big].
\end{equation*}
This is a complex of abelian groups and the naive Chow group of $\Gamma$ is defined to be
\begin{equation*}
\CH_{naif}(\Gamma)=H^0({\cal Z}_{\mathrm{mot}}(\Gamma,*)).
\end{equation*}
By construction one sees that for $\Gamma=\bbZ_{X,W}(n)[2n]$ we have an identification with the usual Bloch's cycle complex with support
\begin{equation*}
{\cal Z}_{\mathrm{mot}}(\Gamma,*)=z^n_W(X,*)
\end{equation*}
that induces an identification on the higher Chow groups
\begin{equation*}
\CH_{naif}(\Gamma[-m])=\CH^n_W(X,m).
\end{equation*}
For a general $\Gamma$  a cycle class map is constructed in \cite[Part I, Ch.II, 2.3.6]{MixMotLev}
\begin{equation}\label{CyclClNaif}
\mathrm{cl}_{naif}(\Gamma):\CH_{naif}(\Gamma)\ra \Hom_{\D^{\mathrm{b}}_{\mathrm{mot}}(F)}(\bb1,\Gamma)
\end{equation}
and proved to be an isomorphism in \cite[Part I, Ch.II, 3.3.10]{MixMotLev}
\footnote{More precisely one has a commutative triangle 
\begin{equation}\label{RmqCycleClLev}
\begin{tabular}{c}
\begin{math}
\xymatrix{{\CH_{naif}(\Gamma)}\ar[rd]^{\mathrm{cl}_{naif}(\Gamma)}\ar[dd] & {} & {}\\
{} & {\Hom_{\D^{\mathrm{b}}_{\mathrm{mot}}(F)}(\bb1,\Gamma)}\ar@{=}[r] & {\Hom_{{\cal{DM}}(F)}(\bb1,\Gamma)}\\
{\mathcal{CH}(\Gamma)}\ar[ru]_{\mathrm{cl}(\Gamma)} & {} & {}}
\end{math}
\end{tabular}
\end{equation}
where $\mathcal{CH}(\Gamma)$ is the Chow group of $\Gamma$ defined \cite[Part I, Ch. II, 2.5.2]{MixMotLev} through the Zariski hyperresolutions of $\Gamma$ and $\mathrm{cl}(\Gamma)$ is the cycle class map defined in \cite[Part I, Ch. II, 2.5.4]{MixMotLev}. Over a field one knows by the Zariski descent property of Bloch's cycle complex \cite{BlochB,LevRev} that the vertical arrow in (\ref{RmqCycleClLev}) is an isomorphism and in \cite[Part I, Ch.II, 3.3.10]{MixMotLev} it is proved that the map $\mathrm{cl}(\Gamma)$ is an isomorphism.}. Thus one  has a cycle class isomorphism
\begin{equation}\label{LevCycCl}
\begin{tabular}{c}
\begin{math}
\xymatrix{{\CH^n_W(X,m)}\ar@{=}[rr]\ar[d]^{\mathrm{cl}^{n,m}_{X,W}} & {} & {\CH_{naif}(\Gamma)}\ar[d]^{\mathrm{cl}_{naif}(\Gamma)}\\
{H^{2n-m}(X,\bbZ(n))_{\mathsf{L}}}\ar@{=}[r] & {\Hom_{{\cal{DM}}(F)}(\bb1,\Gamma)}\ar@{=}[r] & {\Hom_{\D^{\mathrm{b}}_{\mathrm{mot}}(F)}(\bb1,\Gamma)}}
\end{math}
\end{tabular}
\end{equation}
when we take $\Gamma=\bbZ_{X,W}(n)[2n-m]$. 

\begin{lemma}\label{LemDescCyclLev}
Let $X$ be a smooth quasi-projective scheme and $\alpha$ be an element in ${\cal Z}^n_W(X)$. We have the equality
\begin{equation*}
\mathrm{cl}_{X,W}^n(\alpha)=[\alpha]_W\circ[1]^{-1}
\end{equation*}
in the cohomology group $H^{2n}_W(X,\bbZ(n))_{\mathsf{L}}$.
\end{lemma}

\begin{proof}
To see that this formula holds one has just to have an insight into the construction of the cycle class map (\ref{CyclClNaif}). Let $\Gamma$ be an object in $\bC^{\mathrm{b}}_{\mathrm{mot}}(F)^*$. The functor (\ref{Zmot}) induces a morphism
\begin{equation}\label{Morcolimcyc}
\begin{tabular}{c}
\begin{math}
\xymatrix@R=.4cm@C=.1cm{{\colim_a\Hom_{\Kho^{\mathrm{b}}_{\mathrm{mot}}(F)}(\fk{e}^{a}\otimes\bb1,\Gamma)}\ar[r] & {\Hom_{\Kho(\bbZ)}(\mathcal{Z}_{\mathrm{mot}}(\fk{e}\otimes\bb1),\mathcal{Z}_{\mathrm{mot}}(\Gamma))}\ar@{=}[d]\\
{} & {\Hom_{\Kho(\bbZ)}(\bbZ,\mathcal{Z}_{\mathrm{mot}}(\Gamma))=H^0(\mathcal{Z}_{\mathrm{mot}}(\Gamma))}}
\end{math}
\end{tabular}
\end{equation}
where the transition morphisms in the colimit are given by the morphisms
\begin{equation*}
\begin{tabular}{c}
\begin{math}
\xymatrix{{\fk{e}^{\otimes b}\otimes\bb1}\ar[r]\ar@{=}[d] & {\fk{e}^{\otimes a}\otimes\bb1}\\
{\fk{e}^{\otimes a}\otimes\fk{e}^{\otimes b-a}\otimes\bb1}\ar[r]^{\raisebox{.8em}{$\scriptstyle{\fk{e}^{\otimes a}\otimes[1]^{\otimes b-a}\otimes\bb1}$}} & {\fk{e}^{\otimes a}\otimes\bb1^{b-a}\otimes\bb1}\ar[u]_{\fk{e}^{\otimes a}\otimes\boxtimes_{\bb1,\ldots,\bb1}}}
\end{math}
\end{tabular}
\qquad b>a
\end{equation*}
The morphism (\ref{Morcolimcyc}) is an isomorphism by \cite[Part I, Ch.I, 3.3.5.iv]{MixMotLev} and we may consider its inverse
\begin{equation*}
\mathrm{cyc}_{\Gamma}:H^0(\mathcal{Z}_{\mathrm{mot}}(\Gamma))\ra\colim_a\Hom_{\Kho^{\mathrm{b}}_{\mathrm{mot}}(F)}(\fk{e}^{a}\otimes\bb1,\Gamma).
\end{equation*}
We have a morphism
\begin{equation*}
\Hom_{\Kho^{\mathrm{b}}_{\mathrm{mot}}(F)}(\fk{e}^{a}\otimes\bb1,\Gamma)\ra\Hom_{\D^{\mathrm{b}}_{\mathrm{mot}}(F)}(\fk{e}^{a}\otimes\bb1,\Gamma)\xrightarrow{-\circ\nu_a^{-1}}\Hom_{\D^{\mathrm{b}}_{\mathrm{mot}}(F)}(\bb1,\Gamma)
\end{equation*}
where $\nu_a$ is the composition
\begin{equation*}
\fk{e}^{\otimes a}\otimes\bb1\xrightarrow{[1]^{\otimes a}\otimes\bb1}\bb1^{\otimes a}\otimes\bb1\xrightarrow{\boxtimes_{\bb1,\ldots,\bb1}}\bb1
\end{equation*}
and from this one has a morphism
\begin{equation*}
\xymatrix{{H^0(\mathcal{Z}_{\mathrm{mot}}(\Gamma))}\ar@/^2em/[rr]^{\mathrm{cl}_{\Gamma}}\ar[r]_(.35){\mathrm{cyc}_{\Gamma}} & {\colim_a\Hom_{\Kho^{\mathrm{b}}_{\mathrm{mot}}(F)}(\fk{e}^{a}\otimes\bb1,\Gamma)}\ar[r] & {\Hom_{\D^{\mathrm{b}}_{\mathrm{mot}}(F)}(\bb1,\Gamma).}}
\end{equation*}
Let $\Gamma_N$ be the object in $\bC^{\mathrm{b}}_{\mathrm{mot}}(F)^*$ with component of homological degree $m$
\begin{equation*}
\Gamma_{N,m}=\bigoplus_{s\in\mathsf{\Delta}^n_{m,\mathrm{nd}}}\Gamma\times\bbZ_{\Delta^m}(0)_{\delta^m}\qquad m\geq 0
\end{equation*}
and differential given by 
{\renewcommand{\arraystretch}{.5}
\begin{equation*}
\xymatrix{{\Gamma_{N,m}}\ar[r]^{d_m} & {\Gamma_{N,m-1}}\\
{\Gamma\times\bbZ_{\Delta^m}(0)_{\delta^m}}\ar[u]^{\begin{tabular}{l}$\scriptstyle{\textrm{component indexed}}$\\$\scriptstyle{\textrm{by the non-degenerate}}$\\$\scriptstyle{\textrm{$m$-simplex $s$}}$\end{tabular}}\ar[r] & {\Gamma\times\bbZ_{\Delta^{m-1}}(0)_{\delta^{m-1}}}\ar[u]_{\begin{tabular}{l}$\scriptstyle{\textrm{component indexed}}$\\$\scriptstyle{\textrm{by the non-degenerate}}$\\$\scriptstyle{\textrm{$m$-simplex $s\circ\delta^{m-1}_i$}}$.\end{tabular}}}
\end{equation*}}\noindent
By \cite[Part I, Ch. II, 2.2.6, 2.2.7]{MixMotLev} the natural map $\Gamma\ra\Gamma\times\bbZ_{\Delta^*}(0)_{\delta^*}$ has a natural factorisation into maps 
\begin{equation*}
\Gamma\xrightarrow{i_N}\Gamma_N\xrightarrow{\Pi_N}\Gamma\times\bbZ_{\Delta^*}(0)_{\delta^*}.
\end{equation*}
The cycle class map (\ref{CyclClNaif}) is then gotten from the morphism $\mathrm{cl}_{\Gamma}$ through the following diagram
\begin{equation*}
\xymatrix{{} & {H^0(\mathcal{Z}_{\mathrm{mot}}(\Gamma))}\ar[d]_{\mathrm{cyc}_{\Gamma}}\ar@/_2em/[ldd]_{i_N\circ-}\ar@/^11em/[dddd]_{\mathrm{cl}_{\Gamma}}\\
{} & {\colim_a\Hom_{\Kho^{\mathrm{b}}_{\mathrm{mot}}(F)}(\fk{e}^{a}\otimes\bb1,\Gamma)}\ar[d]^{i_N\circ-}\\
{\colim_NH^0(\mathcal{Z}_{\mathrm{mot}}(\Gamma_N))}\ar[r]^{\raisebox{1em}{$\scriptstyle{\colim_N\mathrm{cyc}_{\Gamma_N}}$}}\ar[d]_{\Pi_N\circ-}^{\mathrm{iso.}} & {\colim_{N,a}\Hom_{\Kho^{\mathrm{b}}_{\mathrm{mot}}(F)}(\fk{e}^{a}\otimes\bb1,\Gamma_N)}\ar[d]\\
{H^0(\mathcal{Z}_{\mathrm{mot}}(\Gamma,*))}\ar@{=}[d] & {\colim_{N,a}\Hom_{\D^{\mathrm{b}}_{\mathrm{mot}}(F)}(\fk{e}^{a}\otimes\bb1,\Gamma_N)}\ar[d]_{i_N^{-1}\circ-\circ\nu_a^{-1}}\\
{\CH_{naif}(\Gamma)}\ar[r]^{\mathrm{cl}_{naif}(\Gamma)} & {\Hom_{\D^{\mathrm{b}}_{\mathrm{mot}}(F)}(\bb1,\Gamma).}}
\end{equation*}
Therefore one has only to check that $\mathrm{cl}_{\Gamma}(\alpha)=[\alpha]_W\circ[1]$ for $\Gamma=\bbZ_{X,W}(n)[2n]$. Letting $j:U\hookrightarrow X$ be the open complement of $W$ one has 
\begin{equation*}
\Gamma=\mathsf{Mc}(\Gamma'\ra\Gamma'')
\end{equation*}
where $\Gamma'=\bbZ_X(n)[2n]$ and $\Gamma''=\bbZ_U(n)[2n]$. By \cite[Part I, Ch. I, 3.2.2.iv]{MixMotLev} the morphism induced by the functor $\mathcal{Z}_{\mathrm{mot}}$
\begin{equation*}
\Hom_{\Kho^{\mathrm{b}}_{\mathrm{mot}}(F)}(\fk{e}\otimes\bb1,\Gamma)\ra H^0(\mathcal{Z}_{\mathrm{mot}}(\Gamma))
\end{equation*}
sends the map 
\begin{equation*}
\begin{pmatrix}\boxtimes_{\Gamma',\bb1}\circ\big([\alpha]\otimes\bb1\big)\\\boxtimes_{\Gamma'',\bb1}\circ\big(h_{X,U,[\alpha],j^*}\otimes\bb1\big)
\end{pmatrix}
\end{equation*}
to the element given by
\begin{equation*}
\begin{pmatrix}
\alpha\\ 0
\end{pmatrix}.
\end{equation*}
This implies that
\begin{equation*}
\mathrm{cyc}_{\Gamma}(\alpha)=\boxtimes_{\Gamma,\bb1}\circ([\alpha]_W\otimes\bb1)
\end{equation*}
and to conclude that $\mathrm{cl}_{\Gamma}(\alpha)=[\alpha]_W\circ[1]$ one has just to remark that the following diagram is commutative in $\D^{\mathrm{b}}_{\mathrm{mot}}(F)$
\begin{equation*}
\xymatrix{{\bb1\otimes\bb1}\ar[r]^{\boxtimes_{\bb1,\bb1}} & {\bb1}\ar[rd]_(.3){[\alpha]_W\circ[1]^{-1}} & {\fk{e}}\ar[l]_{[1]}\ar[d]^{[\alpha]_W}\\
{\fk{e}\otimes\bb1}\ar[r]_{[\alpha]_W\otimes\bb1}\ar[u]^{[1]\otimes\bb1} & {\Gamma\otimes\bb1}\ar[r]_(.6){\boxtimes_{\Gamma,\bb1}} & {\Gamma.}}
\end{equation*}
\end{proof}

\noindent We use the notation of \cite[Part I, Ch.VI, 1.1.1]{MixMotLev}. Thus $\square^1$ is the open subscheme of $\mathbb{P}^1$ complement of the unit section and $\square^n$ denotes the product of $n$ copies of $\square^1$. Let $(t_1,\ldots, t_n)$ be the natural coordinates on $\square^n$, we denote by $\partial\square^n$ the collection of divisors 
$$\partial\square^n=\{t_1=0,t_1=\infty,\ldots,t_n=0,t_n=\infty\}.$$ 

\begin{proposition}\label{CompMotCoh}
Let $X$ be a smooth quasi-projective scheme and $n,m$ be two integers with $n$ non negative. The following diagram is commutative
\begin{equation*}
\xymatrix{{\CH^{n}(X,m)}\ar[r]^{\mathrm{cl}^{n,m}_X}\ar[rd]_{\fk{cl}_X^{n,m}} & {H^{2n-m}(X,\bbZ(n))_{\mathsf{L}}}\ar[d]^{\Upsilon_{\fksh}^{\mathrm{eff}}}\\
{} & {H^{2n-m}(X,\bbZ(n))_{\mathsf{V}}.}}
\end{equation*}
\end{proposition}

\begin{proof}
We proceed as in the proof of the proposition 3.8 of \cite{IvoB}. Let $W$ be a closed subscheme of $X$. Remark first that one has a commutative triangle
\begin{equation*}
\xymatrix{{\CH^{n}_W(X)}\ar[r]^{ \mathrm{cl}^{n}_{X,W}}\ar[rd]_{ \fk{cl}_{X,W}^{n}} & { H^{2n}_W(X,\bbZ(n))_{\mathsf{L}}}\ar[d]\\
{} & { H^{2n}_W(X,\bbZ(n))_{\mathsf{V}}.}} 
\end{equation*}
Indeed as the map $\mathcal{Z}^n_W(X)\ra\CH^n_W(X)$ is surjective, one has only to check that for a codimension $n$ cycle $\alpha$ on $X$ with support in $W$ 
\begin{equation*}
\Upsilon^{\mathrm{eff}}_{\fksh}\Big( \mathrm{cl}^{n}_{X,W}(\alpha)\Big)= \fk{cl}^{n}_{X,W}(\alpha)
\end{equation*}
and this equality is a simple consequence of lemma \ref{LemDescCyclLev} since 
\begin{equation*}
\begin{split}
\Upsilon^{\mathrm{eff}}_{\fksh}\Big( \mathrm{cl}^{n}_{X,W}(\alpha)\Big)&=\Upsilon^{\mathrm{eff}}_{\fksh}\Big( [\alpha]_W\circ [1]^{-1}\Big)=\Upsilon^{\mathrm{eff}}_{\fksh}\Big( [\alpha]_W\Big)\circ\Upsilon^{\mathrm{eff}}_{\fksh}\Big( [1]^{-1}\Big)\\
&=\Upsilon^{\mathrm{eff}}_{\fksh}\Big( [\alpha]_W\Big)= \fk{cl}^{n}_{X,W}(\alpha).
\end{split}
\end{equation*} 
One has canonical isomorphisms
\begin{eqnarray*}
\CH^n(X,m) & = & \CH^n(X\square^m;X\partial\square^m)\\
 H^{2n-m}(X,\bbZ(n))_{\mathsf{L}} & = &  H^{2n}(X\square^m;\partial\square^m,\bbZ(n))_{\mathsf{L}}\\
 H^{2n-m}(X,\bbZ(n))_{\mathsf{V}} & = &  H^{2n}(X\square^m;\partial\square^m,\bbZ(n))_{\mathsf{V}}.
\end{eqnarray*}
and so has just to check that the following diagram is commutative
\begin{equation*}
\xymatrix@C=2cm{{\CH^n(X;D)}\ar[r]^{ \mathrm{cl}_{(X;D)}^{n}}\ar[rd]_{ \fk{cl}_{(X;D)}^{n}} & { H^{2n}(X;D,\bbZ(n))_{\mathsf{L}}}\ar[d]^{\Upsilon^{\mathrm{eff}}_{\fksh}}\\
{} & { H^{2n}(X;D,\bbZ(n))_{\mathsf{V}}}}
\end{equation*}
for all subschemes $D_1,\ldots,D_r$ forming a normal crossing divisor in $X$. Fix a pure codimension $n$ closed subscheme $W$ of $X$ which intersects all the $D_i$ properly. Using the description of the relative Chow groups given above lemma 2.5 in \cite{LevRev} one sees that one has only to prove that the following diagram is commutative
\begin{equation}\label{ChowRelSupp}
\begin{tabular}{c}
\begin{math}
\xymatrix@C=2cm{{\CH^n_W(X;D)}\ar[r]^{ \mathrm{cl}_{(X;D),W}^{n}}\ar[rd]_{ \fk{cl}_{(X;D),W}^{n}} & { H^{2n}_W(X;D,\bbZ(n))_{\mathsf{L}}}\ar[d]^{\Upsilon^{\mathrm{eff}}_{\fksh}}\\
{} & { H^{2n}_W(X;D,\bbZ(n))_{\mathsf{V}}.}}
\end{math}
\end{tabular}
\end{equation}
From relativization long exact sequences, their compatibility with cycle class maps and semi-purity we have a diagram whose rows are exact sequences
\begin{equation*}
\xymatrix@C=.4cm{
{0}\ar[r] & { H^{2n}_W(X;D,\bbZ(n))_{\mathsf{L}}}\ar[r] & { H^{2n}_W(X,\bbZ(n))_{\mathsf{L}}}\ar[r] & {\bigoplus_{i=1}^r H^{2n}_{W\cap D_i}(D_i,\bbZ(n))_{\mathsf{L}}}
\\
{0}\ar[r] & {\CH^n_{W}(X;D)}\ar[r]\ar[d]_{ \fk{cl}_{(X;D),W}^{n}}\ar[u]^{ \mathrm{cl}_{(X;D),W}^{n}} & {\CH^n_W(X)}\ar[r]\ar[d]_{ \fk{cl}_{X,W}^{n}}\ar[u]^{ \mathrm{cl}_{X,W}^{n}} & {\bigoplus_{i=1}^r\CH^n_{W\cap D_i}(D_i)}\ar[d]_{\oplus_{i=1}^r \fk{cl}^{n}_{D_i,W\cap D_i}}\ar[u]^{{\oplus_{i=1}^r \mathrm{cl}^{n}_{D_i,W\cap D_i}}}\\
{0}\ar[r] & { H^{2n}_W(X;D,\bbZ(n))_{\mathsf{V}}}\ar[r] & { H^{2n}_W(X,\bbZ(n))_{\mathsf{V}}}\ar[r] & {\bigoplus_{i=1}^r H^{2n}_{W\cap D_i}(D_i,\bbZ(n))_{\mathsf{V}}}}
\end{equation*}
and the commutativity of  (\ref{ChowRelSupp}) follows from the compatibility of the functor $\Upsilon^{\mathrm{eff}}_{\fksh}$ with relativization long exact sequences.
\end{proof}

\begin{corollary}\label{CorPlFid}
Let $\Gamma$ be an object in $\mathcal{DM}^{\mathrm{eff}}_{\fksh}(F)$ and $n$ a non negative integer. The induced map 
\begin{equation*}
\Hom_{\mathcal{DM}^{\mathrm{eff}}_{\fksh}(F)}\big(\bbZ(n),\Gamma\big)\xrightarrow{\Upsilon^{\mathrm{eff}}_{\fksh}}\Hom_{DM^{\mathrm{eff}}_{-}(F)}\big(\bbZ(n), \Upsilon^{\mathrm{eff}}_{\fksh}(\Gamma)\big)
\end{equation*}
is an isomorphism. 
\end{corollary}

\begin{proof}
We may assume that $\Gamma=\bbZ_X(r)[m]$ where $X$ is a smooth quasi-projective scheme and $n,m$ are integers with $n$ non negative. Since 
\begin{equation*}
\Upsilon^{\mathrm{eff}}_{\fksh}\big(\Gamma\big)=\Homi(M(X),\bbZ(r))[m]
\end{equation*}
by Voevodsky's cancellation theorem \cite{Cancel} we have
\begin{equation*}
\begin{split}
\Hom_{DM^{\mathrm{eff}}_{-}(F)}\left(\bbZ(n),\Upsilon^{\mathrm{eff}}_{\fksh}\big(\Gamma\big)\right) &=\Hom_{DM^{\mathrm{eff}}_{-}(F)}\big(M(X)(n),\bbZ(r)[m]\big)\\
&=\Hom_{DM_{-}(F)}\big(M(X)(n),\bbZ(r)[m]\big)\\
&=\Hom_{DM_{-}(F)}\big(M(X),\bbZ(r-n)[m]\big)\\
&=H^m(X,\bbZ(r-n))_{\mathsf{V}}
\end{split}
\end{equation*}
and similarly
\begin{equation*}
\begin{split}
\Hom_{\mathcal{DM}^{\mathrm{eff}}_{\fksh}(F)}\big(\bbZ(n),\Gamma\big) &=\Hom_{\mathcal{DM}_{\fksh}(F)}\big(\bbZ(n),\bbZ_X(r)[m]\big)\\
&=\Hom_{\mathcal{DM}_{\fksh}(F)}\big(\bbZ,\bbZ_X(r-n)[m]\big)\\
&=H^m(X,\bbZ(r-n))_{\mathsf{L}}
\end{split}
\end{equation*}
By the cycle class maps (\ref{VoeCycCl}) and (\ref{LevCycCl}) we may assume $r>n$ since otherwise both are zero. The result follows then from proposition \ref{CompMotCoh}.
\end{proof}

\section{Comparison of triangulated categories}\label{SecComTr}

As mentionned in the introduction, the functor obtained in proposition \ref{Upsilonsheff} is not the right one since it is not a tensor functor and does not extend to non effective triangulated motives. We shall now explain how to get rid of these problems. Consider for that the restriction 
\begin{equation*}
\Upsilon^{\mathrm{eff}+}_{\fksh}:\mathcal{DM}^{\mathrm{eff}+}_{\fksh}(F)\hookrightarrow \mathcal{DM}^{\mathrm{eff}}_{\fksh}(F)\xrightarrow{\Upsilon^{\mathrm{eff}}_{\fksh}}DM^{\mathrm{eff}}_{-}(F). 
\end{equation*}
The fundamental properties of this triangulated commutative external product are given by the following proposition.
\begin{proposition}\label{FundThm}
Assume either that
\begin{itemize}
\item{$\mathrm{char}(F)=0$ and $A=\bbZ$, or}
\item{$\mathrm{char}(F)>0$ and $A=\bbQ$.}
\end{itemize}
The triangulated commutative external product 
\begin{equation*}
\mathcal{DM}^{\mathrm{eff}+}_{\fksh}(F,A)\xrightarrow{\Upsilon^{\mathrm{eff}+}_{\fksh,A}} DM^{\mathrm{eff}}_{-}(F,A). 
\end{equation*}
is a tensor functor which takes its values in $DM^{\mathrm{eff}}_{gm}(F,A)$. In addition it induces a commutative square of triangulated tensor functors
\begin{equation*}
\xymatrix{{\mathcal{DM}^+_{\fksh}(F,A)}\ar[r]^{\Upsilon^+_{\fksh,A}} & {DM_{gm}(F,A)}\\
{\mathcal{DM}^+_{\fksh}(F)^{\mathrm{pr}} }\ar[r]^{\Upsilon^+_{\fksh}}\ar[u] & {DM_{gm}(F)^{\mathrm{pr}}}\ar[u]}
\end{equation*}
where the horizontal arrows are equivalences of categories.
\end{proposition}

\begin{proof}
Recall that the category $\mathcal{DM}^{\mathrm{eff}+}_{\fksh}(F)$ is the pseudo-abelian hull of the full triangulated subcategory $\D^{\mathrm{b},\mathrm{eff}+}_{\mathrm{mot}\fksh}(F)$ of $\D^{\mathrm{b},\mathrm{eff}}_{\mathrm{mot}\fksh}(F)$ generated by the motives $\bbZ_X(n)$ where $X$ is a smooth quasi-projective scheme of dimension at most $n$. By \cite[Corollary B.2]{HuberKahn} for such a motive one has 
\begin{equation*}
\Upsilon^{\mathrm{eff}+}_{\fksh,A}(\bbZ_X(n))=\Homi(M(X),A(n))=M^*(X)(n)
\end{equation*}
which belongs to $DM^{\mathrm{eff}}_{gm}(F,A)$. This proves that $\Upsilon^{\mathrm{eff}+}_{\fksh,A}$ takes its values in this triangulated category. Futhermore by construction this functor is gotten from the functor
\begin{equation*}
\Upsilon^{\mathrm{eff}}_{\mathrm{mot}}\circ i_{\mathrm{mot}}\circ r_{\mathrm{mot}}:\Kho^{\mathrm{b},\mathrm{eff}+}_{\mathrm{mot}\fksh}(F)\ra DM^{\mathrm{eff}}_{-}(F,A).
\end{equation*}
For $\Gamma=\bbZ_X(n)$ and $\Delta=\bbZ_Y(m)$ where $X,Y$ are smooth quasi-projective schemes of dimension at most $n,m$, we have a commutative square
\begin{equation}\label{SqKun}
\xymatrix{{\Upsilon_{\mathrm{mot}\fksh}^{\mathrm{eff}}(\Gamma\otimes\Delta)}\ar[r]^{\Upsilon_{\mathrm{mot}\fksh}^{\mathrm{eff}}(\boxtimes_{\Gamma,\Delta}^{\fksh})} & {\Upsilon_{\mathrm{mot}\fksh}^{\mathrm{eff}}(\Gamma\times\Delta)}\\
{\Upsilon_{\mathrm{mot}\fksh}^{\mathrm{eff}}(\Gamma)\otimes\Upsilon_{\mathrm{mot}\fksh}^{\mathrm{eff}}(\Delta)}\ar@{-}[d]_{\mathrm{iso.}}\ar[u]^{\mathrm{iso.}} & {}\\
{M^*(X)(n)\otimes M^*(Y)(m)}\ar[r]^{(\ref{SqKun})} & {M^*(X\times Y)(n+m).}\ar@{-}[uu]^{\mathrm{iso.}}}
\nonumber
\end{equation}\addtocounter{equation}{1}\noindent
where the map (\ref{SqKun}) is an isomorphism since $DM_{gm}(F,A)$ is a rigid tensor category. Remark \ref{CondTens} shows that $\Upsilon^{\mathrm{eff}+}_{\fksh,A}$ is a tensor functor. Now since for a non negative integer $n$
\begin{equation*}
\Upsilon^{\mathrm{eff}+}_{\fksh,A}(\bbZ(n))=\bbZ(n)
\end{equation*}
we see that $\Upsilon^{\mathrm{eff}+}_{\fksh,A}$ induce a triangulated tensor functor
\begin{equation*}
\Upsilon^{+}_{\fksh,A}:\mathcal{DM}^{+}_{\fksh}(F)\ra DM_{gm}(F,A)
\end{equation*}
It remains to check that this functor is an equivalence. Corollary \ref{CorPlFid} implies that the map 
\begin{equation*}
\Hom_{\mathcal{DM}_{\fksh}(F,A)}(\bb1,\Gamma)\xrightarrow{\Upsilon^+_{\fksh,A}}\Hom_{DM_{gm}(F,A)}(\bbZ,\Upsilon_{\fksh,A}(\Gamma))
\end{equation*}
is an isomorphism for any object $\Gamma$ in $\mathcal{DM}^{+}_{\fksh}(F)$. By \cite[Part I, Ch. IV, 1.5.1]{MixMotLev} this implies that $\Upsilon_{\fksh,A}$ is fully faithful. Let $X$ be a smooth quasi-projective scheme and $n$ an integer. For any integer $d$ greather than the dimension of $X$ the motive $\bbZ_X(d)$ belongs to $\mathcal{DM}^{\mathrm{eff}+}_{\fksh}(F)$ and thus the isomorphism
\begin{equation*}
\begin{split}
\Upsilon^+_{\fksh,A}\big(\bbZ_X(d)(n-d)\big)& =\Upsilon^{\mathrm{eff}+}_{\fksh,A}\big(\bbZ_X(d)\big)(n-d)=\big(M^*(X)(d)\big)(n-d)\\
&=M^*(X)(n)
\end{split}
\end{equation*}
shows that the essential image is a pseudo-abelian triangulated subcategory that contains the $M^*(X)(n)$ and so is all of $DM_{gm}(F,A)$. With a similar proof one checks that $\Upsilon^+_{\fksh,A}$ induces a triangulated tensor equivalence
\begin{equation*}
\mathcal{DM}^+_{\fksh}(F)^{\mathrm{pr}}\xrightarrow{\Upsilon^+_{\fksh}}DM_{gm}(F)^{\mathrm{pr}}.
\end{equation*}
This proves the proposition.
\end{proof}

\noindent Now choose a tensor triangulated quasi-inverse $\rho$ to the triangulated tensor functor $\mathcal{DM}^{+}_{\fksh}(F)\ra\mathcal{DM}(F) $
and denote by $\Upsilon$ the composition
\begin{equation*}
\Upsilon:\mathcal{DM}(F,A)\xrightarrow{\rho_A}\mathcal{DM}^{+}_{\fksh}(F,A)\xrightarrow{\Upsilon^{+}_{\fksh,A}} DM_{gm}(F,A).
\end{equation*}
Then proposition \ref{FundThm} has the following consequence which is the main result of the paper:
\begin{theorem}\label{MainTheorem}
Assume either that
\begin{itemize}
\item{$\mathrm{char}(F)=0$ and $A=\bbZ$, or}
\item{$\mathrm{char}(F)>0$ and $A=\bbQ$.}
\end{itemize}
Then  the triangulated tensor functor $$\mathcal{DM}(F,A)\xrightarrow{\Upsilon}DM_{gm}(F,A)$$ is an equivalence such that for a quasi-projective scheme $X$ and an integer $n$ we have an isomorphism 
\begin{equation*}
\Upsilon(\bbZ_X(n))\simeq M(X)^*(n).
\end{equation*}
In addition the equivalence $\Upsilon$ induces a triangulated tensor equivalence
\begin{equation*}
{\cal DM}(F)^{\mathrm{pr}}\xrightarrow{\Upsilon} DM_{gm}(F)^{\mathrm{pr}}.
\end{equation*}
\end{theorem}

\appendix

\section*{Appendix}
\setcounter{section}{1}

\subsection{Cosimplicial and simplicial objects}\label{AppendixDelta}

We denote the category of standard simplexes by $\mathsf{\Delta}$. It is the category with objects the ordered sets $\mathsf{\Delta}^n=\{0<\cdots <n\}$ and with morphisms order-preserving maps. The map in $\mathsf{\Delta}$ are generated by the coface maps
\begin{equation*}
\delta^n_k:\mathsf{\Delta}^n\ra\mathsf{\Delta}^{n+1}\qquad\delta^n_k(i)= \begin{cases} i & \textrm{if $i\in\{0,\ldots,k-1\}$}\\ i+1 & \textrm{if $i\in\{k,\ldots,n\}$}  \end{cases}
\end{equation*}
indexed by $k=0,\ldots,n+1$ and the codegeneracy maps
\begin{equation*}
\sigma^n_k:\mathsf{\Delta}^n\ra\mathsf{\Delta}^{n-1}\qquad\sigma^n_k(i)= \begin{cases} i & \textrm{if $i\in\{0,\ldots,k\}$}\\ i-1 & \textrm{if $i\in\{k+1,\ldots,n\}$}  \end{cases}
\end{equation*}
indexed by $k=0,\ldots,n$. In the set of $m$-simplexes of $\mathsf{\Delta}^n$
\begin{equation*}
\mathsf{\Delta}_m^n=\Hom_{\mathsf{\Delta}}(\mathsf{\Delta}^m,\mathsf{\Delta}^n)
\end{equation*}
one denotes by $\mathsf{\Delta}_{m,\mathrm{nd}}^n$ the subset of non-degenerate $m$-simplexes formed by injective maps. Let $p,q$ be non negative integers. Recall that a $(p,q)$-shuffle $s$ is a permutation in $\mathfrak{S}_{p+q}$ such that
\begin{equation*}
s(1) < \cdots < s(p)\qquad\textrm{and}\qquad s(p+1)< \cdots <s(p+q).
\end{equation*}
Such a permutation defines two maps $s^{p,q}_p\in\mathsf{\Delta}^p_{p+q} $ and $s^{p,q}_q\in\mathsf{\Delta}_{p+q}^q$ by
\begin{eqnarray*}
s^{p,q}_p(i) & = & |k\in\{1,\ldots,p\}:k\leq s(i)|\\
s^{p,q}_q(i) & = & |k\in\{p+1,\ldots,p+q\}:k\leq s(i)|.
\end{eqnarray*}
In addition one has two maps $t^{p,q}_p\in\mathsf{\Delta}^{p+q}_p$ and $t^{p,q}_q\in\mathsf{\Delta}^{p+q}_q$ given by 
\begin{eqnarray*}
t^{p,q}_p(i) & = & i\\
t^{p,q}_q(i) & = & p+i.
\end{eqnarray*} 
Let $\mathsf{A}$ be a tensor category. For a cosimplicial object $A\in\mathsf{\Delta}\mathsf{A}$ we denote by $\mathrm{cc}A$ the associated cochain complex with the usual differential
\begin{equation*}
d^n=\sum_{k=0}^{n+1}(-1)^kA(\delta^n_k):A^n\ra A^{n+1}.
\end{equation*}
For cosimplicial objects $A,B\in\mathsf{\Delta}\mathsf{A}$ the \textit{Alexander-Whitney} map
\begin{equation*}
\boxtimes^{\mathrm{cc}}_{A,B}:(\mathrm{cc}A)\otimes(\mathrm{cc}B)\ra\mathrm{cc}(A\otimes B)
\end{equation*}
is the associative but only \textit{commutative up to homotopy} map defined in degree $n$ by
\begin{equation*}
\xymatrix@C=2cm{{\bigoplus_{p+q=n}A^p\otimes B^q}\ar[r]^{\boxtimes^{\mathrm{cc}}_{A,B}} & {A^n\otimes B^n}\\
{A^p\otimes B^q.}\ar[ru]_{\kern2em A(t^{p,q}_p)\otimes B(t^{p,q}_q)}\ar[u] & {}}
\end{equation*}
For a simplicial object $A\in\mathsf{\Delta}^{\op}\mathsf{A}$ we denote by $\mathrm{c}A$ the associated chain complex with the usual differential
\begin{equation*}
d_n=\sum_{k=0}^{n}(-1)^kA(\delta^{n-1}_k):A_n\ra A_{n-1}.
\end{equation*}
For simplicial objects $A,B\in\mathsf{\Delta}^{\op}\mathsf{A}$ the \textit{shuffle} map
\begin{equation*}
\boxtimes^{\mathrm{c}}_{A,B}:(\mathrm{c}A)\otimes(\mathrm{c}B)\ra\mathrm{c}(A\otimes B)
\end{equation*}
is the associative and \textit{commutative} map \cite[theorem 5.2]{EilenMcLane} defined in degree $n$ by
\begin{equation*}
\xymatrix@C=2cm{{\bigoplus_{p+q=n}A^p\otimes B^q}\ar[r]^{\boxtimes^{\mathrm{c}}_{A,B}} & {A^n\otimes B^n}\\
{A^p\otimes B^q}\ar[ru]_{\kern4em \sum_s\epsilon(s)\left[A(s^{p,q}_p)\otimes B(s^{p,q}_q)\right]}\ar[u] & {}}
\end{equation*}
where the sum is taken over all $(p,q)$-shuffles $s$.

\subsection{Nisnevich sheaves with transfers}\label{AppendixNtr}

Denotes by $\Ptr$ the category of presheaves with transfers on $\Sm_F$. The tensor product on $\Ptr $ is obtained by the formula
\begin{equation*}
{\scr F}\otimesptr{\scr G}=\colim_{X\in\Sm_F/{\scr F}}\colim_{Y\in\Sm_F/{\scr G}}\bbZtr[X\times Y].
\end{equation*}
and the one on $\Ntr $ is obtained by taking the associated Nisnevich sheaf and denoted $\otimestr$. Recall \cite[\S2]{SusVoeA} that for preasheaves with transfers $\scr F,\scr G$ and a quasi-projective scheme $X$ one has
\begin{equation*}
(\scr F\otimesptr \scr G)(X)=\left(\bigoplus_{Y,Z\in\Sm_F}{\scr F}(Y)\otimes{\scr G}(Z)\otimes\Corrf(X,Y\times Z)\right)/R
\end{equation*}
where $R$ is the subgroup by the elements
\begin{equation*}
\alpha^*s_{\mtiny{Y}}\otimes \beta^*s_{\mtiny{Z}}\otimes \gamma-s_{\mtiny{Y}}\otimes s_{\mtiny{Z}}\otimes(\alpha\otimes\beta)\circ \gamma
\end{equation*}
with $\alpha\in\Corrf(Y',Y)$, $\beta\in\Corrf(Z',Z)$, $\gamma\in\Corrf_S(X,Y'\times_SZ')$, $s_{\mtiny{Y}}\in{\scr F}(Y)$ and $s_{\mtiny{Z}}\in{\scr G}(Z)$.
The internal Hom in $\Ntr$ is given by the formula
\begin{equation*}
\Homotr(\scr F,\scr G)(X)=\Hom({\scr F}\otimes\bbZtr[X],\scr G)
\end{equation*}
for $X$ a quasi-projective scheme.

\par\bigskip\noindent 
{\bf{Acknowledgements.}} The whole paper owes a lot to Marc Levine's work of motives and it is a pleasure to thank him for its useful insights.
The author has benefitted during this research from support of the European Commission through its FP6 program and its Marie Curie Intra-European Fellowship contract No. 041076.

\vspace{2cm}
\begin{flushright}
\begin{tabular}{l}
	Florian Ivorra\\
	FB Mathematik\\
	Universit\"at Duisburg-Essen\\
	45117 Essen\\
	Germany\\
	florian.ivorra@uni-due.de
\end{tabular}
\end{flushright}
\end{document}